\RequirePackage{fix-cm}
\documentclass[smallextended]{svjour3}       % onecolumn (second format)
\smartqed  % flush right qed marks, e.g. at end of proof
\usepackage{cite}
\usepackage{graphicx}
\usepackage{graphicx,amssymb,amsmath}
\usepackage{booktabs}
\usepackage{array}
\usepackage{color}
\usepackage{enumerate}
\usepackage{algorithm}
\usepackage{algorithmic}
\usepackage{float}
\usepackage{multirow}
\usepackage{rotating}
\usepackage[misc]{ifsym}
\usepackage{bm}
\usepackage[colorlinks,linkcolor=blue,citecolor=black]{hyperref}
\usepackage{amsmath}
\begin{document}

\title{APP-Hom Method for Box Constrained Quadratic Programming%\thanks{Grants or other notes
}
\subtitle{}


\author{Guoqiang Wang \and Bo Yu\and Zixuan Chen}


\institute{Guoqiang Wang \at
School of Mathematical Sciences, Dalian University of Technology, Dalian, Liaoning 116024, P.
R. China             \\
              \email{wangguojim@mail.dlut.edu.cn}           %  \\
           \and
           Bo Yu({\Letter}) \at
School of Mathematical Sciences, Dalian University of Technology, Dalian, Liaoning 116024, P.
R. China           \\
           \email{yubo@dlut.edu.cn}
            \and
           Zixuan Chen \at
School of Mathematical Sciences, Dalian University of Technology, Dalian, Liaoning 116024, P.
R. China           \\
           \email{chenzixuan@mail.dlut.edu.cn}
}

\date{Received: date / Accepted: date}

\maketitle

\begin{abstract}
In this paper, based on a $Q$-linear convergence analysis and an estimate of the linear convergence factor of the proximal point (PP) algorithm  for solving box constrained quadratic programming (BQP) problems, an accelerated proximal point (APP) algorithm for solving BQP problems is presented. To solve the strictly convex BQP problems in each step of the APP algorithm, an efficient homotopy method, which tracks the solution path of a parametric quadratic program, is given. The algorithm with APP algorithm as outer iteration  and the homotopy method as inner iteration is named by APP-Hom. The inner homotopy method is efficient by implementing, a warm-start technique based on the accelerated proximal gradient (APG) method, an $\varepsilon$-relaxation technique for checking prime and dual feasibility and determining/correcting the active set. Numerical tests for randomly generated dense and sparse BQPs, BQPs arising from image deblurring, BQPs in SVM, as well as discretized obstacle problem, elastic-plastic torsion problem, and the journal bearing problem show that the APP algorithm takes much less steps than the PP algorithm, the homotopy method is very efficient for strictly-convex BQP, and in consequence, that the APP-Hom is very efficient for non-convex BQP.
\keywords{ box constrained quadratic programming\and non-convex\and accelerated proximal point algorithm\and homotopy\and support vector machine}
\subclass{90C20 \and 90C26 }
\end{abstract}

\section{Introduction}

In this paper, we consider box constrained quadratic program (BQP)
\begin{equation}\label{equation1.1}
\begin{array}{l}
\min ~~ q(x) \\	
\rm{s.~t.}~~~\emph{l} \leq \emph{x}\leq \emph{u},
\end{array}
\end{equation}
where $q(x)=\frac{1}{2}x^{T}Qx+r^{T}x$, $Q\in \mathbb{R}^{n\times n}$ is symmetric but maybe indefinite, $r$, $l$ and $u$ are $n$-dimensional column vectors.

A basic approach solving $(\ref{equation1.1})$ is the projected-gradient method \cite{rosen1960gradient}, which was proposed for general optimization with simple constraints. However, this method is slow at the end of the iterations, which makes it unappealing when medium to high accuracy solutions are desired. In such cases, utilizing the second order information might be more preferable, e.g. projected-Newton method \cite{bertsekas1982projected}. Though projected-Newton method exhibits strong convergence properties, computing the second order information is very expensive. So approximate second order information are adopted, e.g., LBFGS-B \cite{byrd1995limited}, and projected-quasi-Newton (PQN-BFGS, PQN-LBFGS) \cite{kim2010tackling}.  Among these projection methods, LBFGS-B,  projected-Newton can solve non-convex BQPs, while PQN-LBFGS requires $q(x)$ is strictly convex. In contrast to LBFGS-B and  projected-quasi-Newton, TRON  \cite{lin1999newton} uses trust region to approximate the second-order information, and has  global and superlinear convergence.

Apart from TRON which retains the constrains, Coleman and Li proposed another kind of trust-region method, called reflective Newton method \cite{coleman1996reflective} which transforms $(\ref{equation1.1})$ to an unconstrained piece quadratic minimization problem by a reflective transformation technique and solves the unconstrained minimization by trust-region Newton method. Reflective Newton method exhibits the similar advantages as TRON: strong convergence properties, global and quadratic convergence; and the same computational costs: solving the trust-region subproblems.

Another basic method for $(\ref{equation1.1})$  is active-set method which solves a sequence of subproblems of the form
\begin{eqnarray}\label{equation1.2}
\min \{q(x^{k}+d)~|d_{j}=0,j\in \mathcal{A}(x^{k})\}, 	
\end{eqnarray}
where $\mathcal{A}(x)=\{j|x_{j}=l_j~or~x_{j}=u_j\}$. Typical active-set algorithm restricts the change in $\mathcal{A}(x^{k})$ and updates the active set by dropping or adding only one constraint along the descent direction at each iteration. Recently, Hager et al. presented a new   active-set algorithm (ASA) \cite{hager2006new} which consists of a nonmonotone gradient projection step, an unconstrained optimization step, and a set of rules for branching between the two steps. ASA is shown to be faster than TRON \cite{lin1999newton} for solving the 50 box constrained problems in the CUTEr library \cite{bongartz1995cute} and competitive with TRON for the 23 box constrained problems in the MINPACK-2 library \cite{averick1992minpack}.

In contrast to active-set method, parametric active-set (PAS) method is essentially a homotopy like method, which was proposed by Ritter \cite{ritter1967method, ritter1981parametric} and Best \cite{best1982algorithm, best1996algorithm}, and implemented in qpOASES \cite{ferreau2014qpoases}. PAS solves a general convex quadratic program
\begin{equation}\label{gqp}
\begin{array}{l}
\min ~\{r^{T}x+\frac{1}{2}x^{T}Q{x} |~l\leq Ax\leq u\},
\end{array}
\end{equation}
which includes (\ref{equation1.1}) as a special case with $A=I$, by tracking the piecewise linear solution path of the following parametric quadratic program (PQP) \begin{equation}\label{pqp}
\begin{array}{l}
\min ~\{g(t)^{T}x+\frac{1}{2}x^{T}Q{x} |~l(t)\leq Ax\leq u(t)\},
\end{array}
\end{equation}
from $t=1$ to $t=0$. The PQP in (\ref{pqp}) is constructed such that: (I) when $t=0$, it becomes (\ref{gqp}), that is to say, $g(0)=r, l(0)=l, u(0)=u$; (II) when $t=1$, its solution $x(1)$ as well as the corresponding multipliers $\lambda_l(1)\ge 0$ and $\lambda_u(1)\ge 0$ with $l(1)\le Ax(1)\le u(1)$, $\lambda_l(1)^T (l(1)-Ax(1))=0$, $\lambda_u(1)^T (Ax(1)-u(1))=0$, and $g(1)=-Qx(1)+A^T(\lambda_l(1)-\lambda_u(1))$ are known. At each step of the homotopy tracking, it needs to solve linear systems:
\begin{equation}\label{paslinearsystem}
\left[\begin{array}{ccc}Q & -A_{\mathcal{A}_l}^T&A_{\mathcal{A}_u}^T\\ -A_{\mathcal{A}_l}&0&0\\ A_{\mathcal{A}_u}&0&0\end{array}\right] \left[\begin{array}{c}x(t) \\\lambda_{\mathcal{A}_l}(t)\\\lambda_{\mathcal{A}_u}(t)\end{array}\right]=
\left[\begin{array}{c}-g(t) \\l_{\mathcal{A}_l}(t)\\u_{\mathcal{A}_u}(t)\end{array}\right],
\end{equation}
which derives from the Karush-Kuhn-Tucker (KKT) conditions, where $\mathcal{A}_l$ and $\mathcal{A}_u$ denote the lower active-set and the upper active-set respectively. Similar to the active set method, the efficiency of PAS depends on the difference between the optimal active sets of the start solution and the target solution, and the size of (\ref{paslinearsystem}). Active-set methods change the active set along the descent directions, while PAS changes along the parameter $t$ from $t=1$ to $t=0$. An advantage of PAS is that the times of changing active set is smaller than active-set method. Whereas, PAS suffers from the same disadvantages as the active-set method that it would be inefficient if the start active set is faraway from that of the target solution. So a good warm-start is important to active-set like methods.

Recently, Nesterov's accelerated proximal gradient algorithms (APG) which take small computational cost at each step and converges in a rate $O(\frac{1}{k^2})$ \cite{nesterov2005smooth, nesterov2007gradient} are widely used to solve general convex optimization with simple constraints. Moreover, it is easy to implement, which makes it be a good alternative method for predicting the optimal active set of (\ref{equation1.1}) when it is convex. However, APG algorithm is slow at the end of the iterations and hard to obtain high-precision solution, which hinders it to be an independent algorithm if high-precision solutions are required. So proper termination criteria are required.

In this paper, we present  competitive algorithms for general BQP. First, based on the framework of proximal point (PP)  algorithm  \cite{rockafellar1976monotone}, we solve non-convex and non-strictly convex BQP by solving a sequence of strictly convex BQP problems. In \cite{rockafellar1976monotone}, Rockafellar showed that PP algorithm has linear convergence rate when $Q$ is positive definite. When $Q$ is not positive definite, Luo and Tseng \cite{luo1993error} proved that PP algorithm  $R$-linearly converges to a stationary point. Compared to the results of Luo and Tseng, we give a deeper analysis on the convergence, that we prove that if the limit point satisfies the strict complementary conditions, then PP algorithm $Q$-linearly converges, more precisely, PP algorithm is essentially linear iteration with the iteration vector located in a eigen-subspace which corresponding eigenvalues are smaller than 1. Furthermore, we show that the limit point is a local minimizer with probability 1 if the initial point is randomly given.  Moreover, according to the convergence analysis above, we give an estimation of the linear convergence factor and present an accelerated PP (APP) algorithm.

Moreover, to effectively implement the PP algorithm and the APP algorithm, we present a simplified and improved PAS algorithm, called homotopy algorithm, for solving strictly-convex BQP (PP and APP subproblems). The homotopy algorithm is an simplified and improved PAS method, and is shown to be much faster than PAS algorithm in qpOASES by adopting two techniques. First, we implement APG for a warm-start; secondly, an $\varepsilon$-relaxation technique for checking prime and dual feasibility and determining/correcting the active set is adopted in the homotopy tracking steps, which ensures the stability of the homotopy algorithm.

Finally, the organization of the rest of this paper is as follows. In the second section, the  optimality conditions of the BQP problem are given. In Section 3, we present a deeper  convergence analysis of PP algortihm. Furthermore, based on the convergence analysis, we present an accelerated PP algorithm. Detail of the algorithms are presented in Section 4. Finally, the numerical experiments are presented in Section 5.

\section{Optimality conditions}
Let
\begin{equation}\label{equation2.1}
\begin{array}{l}
\Omega=\{x\in R^{n}: l\leq x\leq u\},
\end{array}
\end{equation} denote the feasible set of problem $($\ref{equation1.1}$)$ and assume $\inf_{x\in\Omega}q(x)>-\infty$. Moreover, let $X^{*}$ denote the set which consists of all KKT points. Then for any $x^{*}\in X^{*}$, the following optimality conditions hold.
\begin{eqnarray} \label{equation2.2a}
\frac{\partial q}{\partial x_{j}}(x^{*})&\geq&0,~ j\in \mathcal{A}_{l}(x^{*}),     \\\label{equation2.2b}
\frac{\partial q}{\partial x_{j}}(x^{*})&=&0  ,~ j\in \mathcal{A}_{m}(x^{*}),        \\\label{equation2.2c}
\frac{\partial q}{\partial x_{j}}(x^{*})&\leq&0 ,~ j\in \mathcal{A}_{u}(x^{*}),
\end{eqnarray}
where $\mathcal{A}_{l}(x)=\{j|x_{j}=l_{j}\}$, $\mathcal{A}_{m}(x)=\{j|l_j<x_{j}<u_{j}\}$ and $\mathcal{A}_{u}(x)=\{j|x_{j}=u_{j}\}$.

\section{Proximal point algorithm and convergence}
PP algorithm solves problem $(\ref{equation1.1})$ when $Q$ is not positive definite as follows:
\begin{eqnarray} \label{equation3.1}
x^{k+1}=\arg\min~\{q(x)+\frac{\gamma}{2}\Arrowvert x-x^{k} \Arrowvert^{2}|x\in \Omega \},
\end{eqnarray}
where $\gamma=\delta-\lambda_{\min}(Q)$, $\delta>0$ is a constant.

Luo and Tseng \cite{luo1993error} proved that (\ref{equation3.1}) $R$-linearly converges to a stationary point. Based on this, we prove that if the limit point satisfies the strict complementary conditions, then PP algorithm is essentially linear iteration with the iteration vectors located in a eigen-subspace which corresponding eigenvalues are smaller than 1. Moreover, we show that the limit point is a local minimizer with probability 1 if the initial point is randomly given. Furthermore, based on the convergence analysis, we give an estimation of the linear convergence factor and present an accelerated PP algorithm.

\subsection{$Q$-linearly convergence}
We first introduce a lemma.
\begin{lemma} \label{Matrix eigenvector}
Assume $A\in R^{n\times n}$ is positive definite. If $\Arrowvert A^{k}v\Arrowvert\to 0$, then there exists $\rho\in(0,1)$, such that $\Arrowvert A^{k}v\Arrowvert\leq \rho\Arrowvert A^{k-1}v\Arrowvert$. Moreover, $v$ belongs to the eigen-subspace which corresponding eigenvalues are smaller than 1.
\end{lemma}

$Proof.$ Let $\lambda_{1}\geq\lambda_{2}\geq...\geq\lambda_{n}>0$ denote the eigenvalues of $A$ and $v_{1},v_{2},...,v_{n}$ are corresponding unit eigenvectors which satisfy $v_{i}^{T}v_{j}=0, \forall i\neq j$. $v$ has linear representation by $v_{1},v_{2},...,v_{n}$ as follows
$$v=\sum_{i=1}^{n}\alpha_{i}v_{i},$$
then
\begin{eqnarray}
A^{k}v=\sum_{i=1}^{n}\alpha_{i}\lambda_{i}^{k}v_{i}\to 0
\end{eqnarray}
which implies $\alpha_{i}=0$, if $\lambda_{i}\geq 1$, that is, $v$ belongs to the eigen-subspace which corresponding eigenvalues are smaller than 1.  Moreover,

\begin{eqnarray} \nonumber
\Arrowvert A^{k}v\Arrowvert&=&\Arrowvert \sum_{\alpha_{i}\neq 0}\alpha_{i}\lambda_{i}^{k}v_{i}\Arrowvert\\\nonumber
&\leq& \max_{\alpha_{i}\neq 0}(\lambda_{i})\Arrowvert \sum_{\alpha_{i}\neq 0}\alpha_{i}\lambda_{i}^{k-1}v_{i}\Arrowvert\\\nonumber
&=& \rho \Arrowvert A^{k-1}v\Arrowvert
\end{eqnarray}
where $\rho=\max_{\alpha_{i}\neq 0}(\lambda_{i})<1$. ~~~\endproof

\begin{theorem} \label{local linear convergence rate}
Assume $\bar{x}$ is the limit point of iterations (\ref{equation3.1}) and the strict complementarity conditions hold at $\bar{x}$. Then there exists $N>0$, such that, the following items hold for all $k\geq N$.
\begin{enumerate}[\rm{(}i)]
\item $\Arrowvert x^{k+1}-x^{k+2} \Arrowvert \leq \bar{\rho} \Arrowvert x^{k}-x^{k+1} \Arrowvert$ and $\Arrowvert x^{k+1}-\bar{x} \Arrowvert \leq  \bar{\rho} \Arrowvert x^{k}-\bar{x} \Arrowvert$.
\item $q(x^{k})-q(\bar{x}) \leq \frac{\Arrowvert Q_{\bar{\mathcal{A}}_{m}\bar{\mathcal{A}}_{m}} \Arrowvert\bar{\rho}^{2(k-N)}}{2(1-\bar{\rho})^{2}} \Arrowvert x^{N}-x^{N+1} \Arrowvert^{2}$,

\end{enumerate}
where $\bar{\rho}=\max_{\lambda_{i}(M)<1}(\lambda_{i}(M))=\frac{\gamma}{\gamma+\min_{\bar{\lambda}_{i}>0}(\bar{\lambda}_{i})}\in(0,1)$, $\bar{\lambda}_{i}$ is the eigenvalue of $Q_{\bar{\mathcal{A}}_{m}\bar{\mathcal{A}}_{m}}$, and  $\bar{\mathcal{A}}_{m}$ denotes $\mathcal{A}_{m}(\bar{x})$, which is similar to $\bar{\mathcal{A}}_{l}$, $\bar{\mathcal{A}}_{u}$.
\end{theorem}

$Proof.$
Since the strict complementarity conditions hold at $\bar{x}$, there exists $N>0$ such that  $\forall k\geq N$
\begin{eqnarray} \label{active_set_equal}
\mathcal{A}_{l}(x^k)=\bar{\mathcal{A}}_{l},~\mathcal{A}_{m}(x^k)=\bar{\mathcal{A}}_{m},~\mathcal{A}_{u}(x^k)=\bar{\mathcal{A}}_{u}.
\end{eqnarray}

Hence, when $k>N$, we have
\begin{eqnarray} \label{frontback1}
Q_{\bar{\mathcal{A}}_{m}\bar{\mathcal{A}}_{m}}x^{k+1}_{\bar{\mathcal{A}}_{m}}+Q_{\bar{\mathcal{A}}_{m}\bar{\mathcal{A}}_{l}}l_{\bar{\mathcal{A}}_{l}}+Q_{\bar{\mathcal{A}}_{m}\bar{\mathcal{A}}_{u}}u_{\bar{\mathcal{A}}_{u}}+r_{\bar{\mathcal{A}}_{m}}+\gamma(x^{k+1}_{\bar{\mathcal{A}}_{m}}-x^{k}_{\bar{\mathcal{A}}_{m}})&=&0,\\\label{frontback2}
Q_{\bar{\mathcal{A}}_{m}\bar{\mathcal{A}}_{m}}x^{k}_{\bar{\mathcal{A}}_{m}}+Q_{\bar{\mathcal{A}}_{m}\bar{\mathcal{A}}_{l}}l_{\bar{\mathcal{A}}_{l}}+Q_{\bar{\mathcal{A}}_{m}\bar{\mathcal{A}}_{u}}u_{\bar{\mathcal{A}}_{u}}+r_{\bar{\mathcal{A}}_{m}}+\gamma(x^{k}_{\bar{\mathcal{A}}_{m}}-x^{k-1}_{\bar{\mathcal{A}}_{m}})&=&0.
\end{eqnarray}

So we arrive at
\begin{eqnarray}
(Q_{\bar{\mathcal{A}}_{m}\bar{\mathcal{A}}_{m}}+\gamma I)(x^{k+1}_{\bar{\mathcal{A}}_{m}}-x^{k}_{\bar{\mathcal{A}}_{m}})=\gamma(x^{k}_{\bar{\mathcal{A}}_{m}}-x^{k-1}_{\bar{\mathcal{A}}_{m}}),\nonumber
\end{eqnarray}
that is,
\begin{eqnarray} \label{poweriteration}
x^{k+1}_{\bar{\mathcal{A}}_{m}}-x^{k}_{\bar{\mathcal{A}}_{m}}=\bar{M}(x^{k}_{\bar{\mathcal{A}}_{m}}-x^{k-1}_{\bar{\mathcal{A}}_{m}}),
\end{eqnarray}
where $\bar{M}=\gamma(Q_{\bar{\mathcal{A}}_{m}\bar{\mathcal{A}}_{m}}+\gamma I)^{-1}$ is positive definite. Similarly, we have
\begin{eqnarray} \label{poweriteration2}
x^{N+1}_{\bar{\mathcal{A}}_{m}}-\bar{x}_{\bar{\mathcal{A}}_{m}}=\bar{M}(x^{N}_{\bar{\mathcal{A}}_{m}}-\bar{x}_{\bar{\mathcal{A}}_{m}}),
\end{eqnarray}
Then we obtain (i) from Lemma \ref{Matrix eigenvector}. Moreover, we have that $x^{k+1}_{\bar{\mathcal{A}}_{m}}-x^{k}_{\bar{\mathcal{A}}_{m}}, k>N$ and $x^{N+1}_{\bar{\mathcal{A}}_{m}}-\bar{x}_{\bar{\mathcal{A}}_{m}}$ belong to the eigen-subspace of $Q_{\bar{\mathcal{A}}_{m}\bar{\mathcal{A}}_{m}}$ which corresponding eigenvalues are smaller than 1.

Furthermore, when $k\geq N$, we obtain
\begin{eqnarray}\label{objective value estimate}\nonumber
q(x^{k})-q(\bar{x})
&=&\frac{1}{2}(x_{\bar{\mathcal{A}}_{m}}^{k}-\bar{x}_{\bar{\mathcal{A}}_{m}})^{T}Q_{\bar{\mathcal{A}}_{m}\bar{\mathcal{A}}_{m}}(x_{\bar{\mathcal{A}}_{m}}^{k}-\bar{x}_{\bar{\mathcal{A}}_{m}})\\\nonumber
                              &\leq&\frac{1}{2}\Arrowvert Q_{\bar{\mathcal{A}}_{m}\bar{\mathcal{A}}_{m}}\Arrowvert \Arrowvert x^{k}-\bar{x} \Arrowvert^{2}\\\nonumber
                              &\leq& \frac{1}{2}\Arrowvert Q_{\bar{\mathcal{A}}_{m}\bar{\mathcal{A}}_{m}}\Arrowvert\sum_{i=k}^{\infty}\Arrowvert x^{i}-x^{i+1} \Arrowvert^{2}\\\nonumber
                              &\leq& \frac{\Arrowvert Q_{\bar{\mathcal{A}}_{m}\bar{\mathcal{A}}_{m}}\Arrowvert}{2(1-\bar{\rho})^{2}}\Arrowvert x^{k}-x^{k+1} \Arrowvert^{2}\\\nonumber
                              &\leq& \frac{\Arrowvert Q_{\bar{\mathcal{A}}_{m}\bar{\mathcal{A}}_{m}} \Arrowvert\bar{\rho}^{2(k-N)}}{2(1-\bar{\rho})^{2}}\Arrowvert x^{N}-x^{N+1} \Arrowvert^{2}.
\end{eqnarray}~~~\endproof

Next, we show that if the strictly complementary conditions hold at $\bar{x}$ and the initial point is randomly given, then it is a local minimizer with probability 1.

Let $\mathcal{A}_{m}^{k}$ denote $\mathcal{A}_{m}(x^k)$, which is similar to $\mathcal{A}_{l}^{k}, \mathcal{A}_{u}^{k}$.

\begin{theorem} \label{theorem pro}
Assume
\begin{enumerate}[\rm{(}i)]
\item $x^{0}$ is randomly given by a random distribution as: $x^{0}_{i}\sim \xi_{i}$ ;
\item $\forall k\geq 1$, $\mathcal{A}_{m}^{k}\cap\mathcal{A}_{m}^{k-1}\neq \varnothing$.
\end{enumerate}
Then there exist  $B^{k,i}\in \mathbb{R}^{n}$,$c^{k}\in \mathbb{R}^{ |\mathcal{A}_{m}^{k}|}$, such that
\begin{eqnarray}\label{ramdomdistribution}
x^{k}_{[\mathcal{A}_{m}^{k}]_{i}}\sim\sum_{j=1}^{n}B^{k,i}_{j}\xi_{j}+c^{k}_{i}, i=1,...,|\mathcal{A}_{m}^{k}|.
\end{eqnarray}
where $[\mathcal{A}_{m}^{k}]_{i}$ is the $i$-th element of $\mathcal{A}_{m}^{k}$.
\end{theorem}

$Proof.$ From (\ref{frontback1}), we have
\begin{eqnarray}\label{solutionform}\nonumber
x^{k}_{\mathcal{A}_{m}^{k}}=M^{k}(x^{k-1}_{ \mathcal{A}_{m}^{k}}+g^k)
\end{eqnarray}
where $M^{k}=\gamma (Q_{ \mathcal{A}_{m}^{k} \mathcal{A}_{m}^{k}}+\gamma I)^{-1}$ is positive definite and $g^{k}=-\frac{r_{ \mathcal{A}_{m}^{k}}-Q_{ \mathcal{A}_{m}^{k} \mathcal{A}_{l}^{k}}l_{ \mathcal{A}_{l}^{k}}-Q_{ \mathcal{A}_{m}^{k} \mathcal{A}_{u}^{k}}u_{ \mathcal{A}_{u}^{k}}}{\gamma}$.

When $k=1$,
\begin{eqnarray}\nonumber
x^{1}_{\mathcal{A}_{m}^{1}}&=&M^{1}(x^{0}_{ \mathcal{A}_{m}^{1}}+g^1)\\\label{x1}
                           &=&M^{1}\sum_{j=1}^{|\mathcal{A}_{m}^{1}|}(x^{0}_{[\mathcal{A}_{m}^{1}]_{j}}e^{1}_j)+M^{1}g^{1}\\\nonumber
                           &=&\sum_{j=1}^{|\mathcal{A}_{m}^{1}|}(x^{0}_{[\mathcal{A}_{m}^{1}]_{j}}M^{1}e^{1}_j)+M^{1}g^{1}
\end{eqnarray}
where  $e^{k}_{j}\in \mathbb{R}^{|\mathcal{A}_{m}^{k}|}$ is the $j$-th coordinate vector.

From (\ref{x1}), we have
$$
B^{1,i}_{j}=\left\{
\begin{array}{lcl}
M^{1}_{ij},      &      & j\in\mathcal{A}_{m}^{k};\\
0,  &      & else
\end{array} \right.
$$
and $c^{1}=M^{1}g^{1}$.

Without loss of generality, assume (\ref{ramdomdistribution}) holds when $k=k_{1}\geq1$. Similar to (\ref{x1})
\begin{eqnarray}\label{solutionform2}
x^{k_{1}+1}_{\mathcal{A}_{m}^{k_{1}+1}}&=&\sum_{j=1}^{|\mathcal{A}_{m}^{k_{1}+1}|}(x^{k_{1}}_{[\mathcal{A}_{m}^{k_{1}+1}]_{j}}M^{k_{1}+1}e^{k_{1}}_j)+M^{k_{1}+1}g^{k_{1}+1},
\end{eqnarray}
that is, \begin{eqnarray}\label{solutionform3}
x^{k_{1}+1}_{[\mathcal{A}_{m}^{k_{1}+1}]_{i}}=\sum_{j=1}^{|\mathcal{A}_{m}^{k_{1}+1}|}(x^{k_{1}}_{[\mathcal{A}_{m}^{k_{1}+1}]_{j}}M_{ij}^{k_{1}+1}+M_{ij}^{k_{1}+1}g_{j}^{k_{1}+1})
\end{eqnarray}
Since $\mathcal{A}_{m}^{k_{1}+1}\cap\mathcal{A}_{m}^{k_{1}}\neq \varnothing$ and $x^{k_{1}}_{[\mathcal{A}_{m}^{k_{1}}]_{i}}\sim\sum_{j=1}^{n}B^{k_{1},i}_{j}\xi_{j}+c^{k_{1}}_{i}$, we have (\ref{ramdomdistribution}) when $k=k_{1}+1$.

Note that, if for some $k$, $\mathcal{A}_{m}^{k}$ is empty, the second assumption is not satisfied. In that case, $x^{k+1}$ would equal to $x^k$. In fact, if $x^k_j=l_j$, then $Q_j^T x^k+r_j+\gamma(x_j^k-x^{k}_j)\geq Q_j^T x^k+r_j+\gamma(x_j^k-x^{k-1}_j)\geq 0 $; if $x^k_j=u_j$, then $Q_j^T x^k+r_j+\gamma(x_j^k-x^{k}_j)\leq Q_j^T x^k+r_j+\gamma(x_j^k-x^{k-1}_j)\leq 0 $, so $x^{k}=\arg\min~\{q(x)+\frac{\gamma}{2}\Arrowvert x-x^{k} \Arrowvert^{2}|x\in \Omega \}$. Moreover, if $\mathcal{A}_{m}^{k}\neq \varnothing$ but $\mathcal{A}_{m}^{k}\cap\mathcal{A}_{m}^{k-1}=\varnothing$, we add small random perturbation to $x_{\mathcal{A}_{m}^{k}}$ and set $k=0$, then it is easy to see that the consequences of Theorem \ref{theorem pro} would also hold.

\begin{theorem} \label{theorem 2}
Assume the assumptions in Theorem \ref{theorem pro} hold. $\bar{x}$ is the limit point of $\{x^{k}\}$, at which the strict complementarity conditions hold. Then $Q_{\bar{\mathcal{A}}_{m}\bar{\mathcal{A}}_{m}}$ is positive semi-definite with probability 1, that is, $\bar{x}$ is a local minimizer.
\end{theorem}

$Proof.$ From Theorem \ref{local linear convergence rate} and \ref{theorem pro}, we have $N$ such that
$$x_{\bar{\mathcal{A}}_{m}}^{N+1}-x_{\bar{\mathcal{A}}_{m}}^{N}=(\bar{M}-I)x^{N}_{\bar{\mathcal{A}}_{m}}+\bar{M}^{N}g^N.$$
Furthermore,
\begin{eqnarray}\label{poweriteration4}\nonumber
\bar{v}_{s}^{T}(x_{\bar{\mathcal{A}}_{m}}^{N+1}-x_{\bar{\mathcal{A}}_{m}}^{N})&=&\bar{v}_{s}^{T}(\bar{M}-I)x^{N}_{\bar{\mathcal{A}}_{m}}+\bar{v}_{s}^{T}\bar{M}g^N\\\nonumber
&=&\bar{v}_{s}^{T}((\frac{\gamma}{\bar{\lambda}_{s}+\gamma}-1)x^{N}_{\bar{\mathcal{A}}_{m}}+\bar{M}g^N),
\end{eqnarray}
where $\bar{v}_{s}$ is the eigenvector of  $Q_{\bar{\mathcal{A}}_{m}\bar{\mathcal{A}}_{m}}$. Since $\bar{M}g^N$ is not a random variable and $x^{N}_{[\bar{\mathcal{A}}_m]_{i}}$ obeys random distribution as (\ref{ramdomdistribution}), we have that
\begin{eqnarray}\label{condition1}
\bar{v}_{s}^{T}(x_{\bar{\mathcal{A}}_{m}}^{N+1}-x_{\bar{\mathcal{A}}_{m}}^{N})\neq0,~\forall\bar{\lambda}_{s}< 0,
\end{eqnarray}
would hold with probability 1.

Moreover, From Theorem \ref{local linear convergence rate}, we know that
\begin{eqnarray} \nonumber
x^{k+1}_{\bar{\mathcal{A}}_{m}}-x^{k}_{\bar{\mathcal{A}}_{m}}=\bar{M}^{(k-N)}(x^{N+1}_{\bar{\mathcal{A}}_{m}}-x^{N}_{\bar{\mathcal{A}}_{m}}) \to 0
\end{eqnarray}
which implies
\begin{eqnarray} \label{con}
\bar{v}_{s}^{T}(x_{\bar{\mathcal{A}}_{m}}^{N+1}-x_{\bar{\mathcal{A}}_{m}}^{N})=0,~\forall\bar{\lambda}_{s}\leq 0,
\end{eqnarray}
by Lemma \ref{Matrix eigenvector}. From (\ref{condition1}) and (\ref{con}). We have that $\bar{\lambda}_{s}\geq 0, \forall s$ would hold with probability 1,
that is, $Q_{\bar{\mathcal{A}}_{m}\bar{\mathcal{A}}_{m}}$ is positive semi-definite.

Next, we will show $\bar{x}$ is a local minimizer.

Assume $d\in R^{n}$ is a feasible direction of $\Omega$ at $\bar{x}$, then
$$d_{\bar{\mathcal{A}}_{0}}\geq 0,~d_{\bar{\mathcal{A}}_{u}}\leq 0.$$

So if there exists $j\in \mathcal{A}_{0}(\bar{x})\cup \mathcal{A}_{u}(\bar{x})$, such that $d_{j}\neq 0$, then
$$q(\bar{x}+\alpha d)=q(\bar{x})+\alpha d^{T}\triangledown q(\bar{x})+o(\alpha^{2}).$$

we have $d^{T}\triangledown q(\bar{x})>0$ for the strict complementarity conditions hold at $\bar{x}$. So $q(\bar{x}+\alpha d)>q(\bar{x})$ when $\alpha>0$ is small enough.

If $d_{\bar{\mathcal{A}}_{0}\cup \bar{\mathcal{A}}_{u}}=0$, then
\begin{eqnarray}\label{feasible direction}\nonumber
q(\bar{x}+\alpha d)-q(\bar{x})&=&\frac{\alpha^{2}}{2}d^{T}Qd+\alpha d^{T}(Q\bar{x}+r)\\\nonumber
                              &=&\frac{\alpha^{2}}{2}d_{\bar{\mathcal{A}}_{m}}^{T}Q_{\bar{\mathcal{A}}_{m}\bar{\mathcal{A}}_{m}}d_{\bar{\mathcal{A}}_{m}}\\\nonumber
                              &\geq&0.
\end{eqnarray}~~~\endproof

\subsection{Accelerated proximal point algorithm}
From Theorem \ref{local linear convergence rate}, we know that PP algorithm  is linear iteration when $(\mathcal{A}_{l}^{k}, \mathcal{A}_{m}^{k}, \mathcal{A}_{u}^{k})$ remains unchanged. Similar to (\ref{poweriteration2}), we have that
\begin{eqnarray} \nonumber
(Q_{\mathcal{A}^k_{m}\mathcal{A}^k_{m}}+\gamma I)(x^{k}_{\mathcal{A}^k_{m}}-x^{k-1}_{\mathcal{A}^k_{m}})=\gamma (x^{k-1}_{\mathcal{A}^k_{m}}-x^{k-2}_{\mathcal{A}^k_{m}}), \forall K_{1}+2\leq k\leq K_{2},
\end{eqnarray}
Let $ \omega_{k}=\frac{\Arrowvert x^{k}-x^{k-1} \Arrowvert }{\Arrowvert x^{k-1}-x^{k-2} \Arrowvert }$. If $x^{k}-x^{k-1}$ and $x^{k}-x^{k-1}$ tend to linear dependent,  $\omega_{k}$ tends to one constant $\bar{\omega}$, we have
\begin{eqnarray} \nonumber
x^{k}-x^{k-1}\approx \bar{\omega}(x^{k-1}-x^{k-2}),
\end{eqnarray}
which implies
\begin{eqnarray} \nonumber
x^{K_2}\approx x^{k-1}-\frac{1-\bar{\omega}^{K_2-k+1}}{1-\bar{\omega}}(x^{k-1}-x^{k}).
\end{eqnarray}
When $K_2>>k$ and $\bar{\omega}<1$ , we have
\begin{eqnarray} \label{xK2}
x^{K_2}\approx x^{k-1}-\frac{1}{1-\bar{\omega}}(x^{k-1}-x^{k}).
\end{eqnarray}

From $(\ref{xK2})$, we see that we do not need to iterate from $k$ to $K_2$ step by step to obtain $x^{K_2}$, instead, we use (\ref{xK2}) to predict the value of $x^{K_2}$. According to this, we present an accelerated PP (APP) iteration as follows
\begin{eqnarray} \label{accelerate PPA}
x^{k+1}=\arg\min~\{q(x)+\frac{\gamma}{2}\Arrowvert x-\frac{x^{k}-\omega_{k}x^{k-1}}{1-\omega_{k}} \Arrowvert^{2}|x\in \Omega \},
\end{eqnarray}
which is obtained from PP algorithm  by replacing $x^k$ with a prediction of $x^{K_2}$.

We do not directly verify whether $(\mathcal{A}_{l}^{k}, \mathcal{A}_{m}^{k}, \mathcal{A}_{u}^{k})$ remain unchanged, instead, we replace PP iterations (\ref{equation3.1}) by APP iterations (\ref{accelerate PPA}) when
\begin{eqnarray} \label{switch condition}
\omega_{k-j+1}<1~and~|\omega_{k-j+1}-\omega_{k-j}|< \varepsilon, j=1,2,...a,
\end{eqnarray}
where $\varepsilon$ and $a$ are pre-given, because we found that (\ref{switch condition}) is a more effective criterion than verifying whether the free variables do not change. In fact, (\ref{switch condition}) does not require  the free variables of $x^k$ are the same, but APP algorithm also exhibits the affect of acceleration.

If (\ref{switch condition}) holds, after one iteration of (\ref{accelerate PPA}), $x^{k+1}$ quickly tends to $x^{K_2}$, but $x^{k+1}$ may not satisfy (\ref{switch condition}). Thus, we need to switch to PP iterations until (\ref{switch condition}) holds again.

Moreover, we have deduced that $Q_{\bar{\mathcal{A}}_{m}\bar{\mathcal{A}}_{m}}$ is positive semi-definite with probability 1, which implies that  $\lambda_{\max}(\gamma(Q_{\bar{\mathcal{A}}_{m}\bar{\mathcal{A}}_{m}}+\gamma I)^{-1})$ is not bigger than 1. So if $\mathcal{A}^k_{m}=\mathcal{A}^{k+1}_{m}=...=\bar{\mathcal{A}}_{m}$, APP iterations converge to  $\bar{x}$ quickly.

\section{Algorithms}
When $Q$ is not positive definite, PP and APP algorithms solve $(\ref{equation1.1})$ by solving a sequence of strictly convex BQP subproblems. So the efficiency of PP and APP algorithms is based on the solving of the strictly convex BQP subproblems. Moreover, in many applications, such as support vector machine, image processing etc, strictly convex BQP is a fundamental problem. So we present an efficient homotopy method for solving strictly convex BQP.

For convenience, we write the strictly convex BQP in a uniform form
\begin{equation}\label{equation4.1}
\begin{array}{l}
\min ~~ \frac{1}{2}z^{T}Hz+f^{T}z \\	
\rm{s.~t.}~~~\emph{l} \leq \emph{z}\leq \emph{u},
\end{array}
\end{equation}
where  $H$ is positive definite. It is easy to see that (\ref{equation3.1}) is a special case of (\ref{equation4.1}) with $H=Q+\gamma I$, $f=r-\gamma x^k$.

\subsection{ Homotopy method for strictly convex BQP}
In this section, we present the homotopy algorithm, which is an simplified implementation of PAS, and is improved with two important techniques: warm-start: using APG to obtain a good prediction of the optimal active set;  an $\varepsilon$-relaxation technique for checking prime and dual feasibility and determining/correcting the active set, which ensures the stability of the homotopy algorithm.

\textbf{$\bullet$ Warm-start.} Since PAS needs a good warm start, we  use APG algorithm to obtain a prediction of the optimal active set.

We know from  $(\ref{equation2.2a})$-$(\ref{equation2.2c})$ that $\bar{z}$ is the solution of $(\ref{equation4.1})$ if and only if $\bar{z}$ satisfies
\begin{eqnarray}\label{equation4.2a}
H_{j}^{T}\bar{z}+f_{j}&\geq&0,~j\in {\mathcal{A}}_{l}(\bar{z}),  \\	\label{equation4.2b}
H_{j}^{T}\bar{z}+f_{j}&=&0, ~j\in {\mathcal{A}}_{m}(\bar{z}), \\	\label{equation4.2c}
H_{j}^{T}\bar{z}+f_{j}&\leq&0, ~j\in {\mathcal{A}}_{u}(\bar{z}),
\end{eqnarray}
where $H_{j}$ denotes the $j$-th column of $H$, $f_{j}$ denotes the $j$-th component of $f$.

Let $\zeta=1$, $\rho_{1}=1$, and $z^{1}=y^{0}\in \mathbb{R}^{n}$ is given, then APG algorithm iterates as follows:
\begin{eqnarray}\label{equation4.22a}
&&y^{\zeta}=\arg\min_{l\leq z\leq u} \langle Hz^{\zeta}+f,z\rangle+\frac{L}{2}\Arrowvert z-z^{\zeta} \Arrowvert^{2}, \\	\label{equation4.22b}
&&\rho_{\zeta+1}=\frac{1+\sqrt{1+4\rho_{\zeta}^2}}{2},\\ \label{equation4.22c}
&&z^{\zeta+1}=y^{\zeta}+(\frac{\rho_{\zeta}-1}{\rho_{\zeta+1}})(y^{\zeta}-y^{\zeta-1}),
\end{eqnarray}
where $L\geq \Arrowvert H \Arrowvert$. At each iteration, $(\ref{equation4.13a})$ can be fast solved by a projection operator as follows:
\begin{equation} \label{equation4.23}
\begin{array}{l}
y^{\zeta}=\mathcal{P}_{[l,u]}(z^{\zeta}-\frac{1}{L}(Hz^{l}+f)).
\end{array}
\end{equation}
So we just need to do one matrix-vector multiplication. Furthermore,
\begin{equation} \label{matrix vector simplicfication}
\begin{array}{l}
Hz=H^T_{\mathcal{A}_{l}(z)}z_{\mathcal{A}_{l}(z)}+H^T_{\mathcal{A}_{m}(z)}z_{\mathcal{A}_{m}(z)}+H^T_{\mathcal{A}_{u}(z)}z_{\mathcal{A}_{u}(z)}.
\end{array}
\end{equation}
A small but important technique should be mentioned that  when $\mathcal{A}_{l}, \mathcal{A}_{u}$ change small, we can compute the first and the third part of (\ref{matrix vector simplicfication}) by just computing the changed indices of $\mathcal{A}_{l}, \mathcal{A}_{u}$. This is very important when the number of the free variables is small.

We terminate APG iterations when  one of the following criteria is satisfied.
\begin{eqnarray} \label{equation4.21a}
&& \pi_{\varepsilon_{1}}(y^{\zeta}) =\pi_{\varepsilon_{1}}(y^{\zeta-i}) ,~for~ i=1,..,S_{\max},\\\label{equation4.21b}
&&\frac{\Arrowvert y^{\zeta}-y^{\zeta-1} \Arrowvert}{\Arrowvert y^{\zeta} \Arrowvert}<\varepsilon_{2},
\end{eqnarray}
where $\pi_{\varepsilon_{1}}(y)=|\{j|l_j+\| y \|\varepsilon_{1}<y_{j}<u_j-\| y \|\varepsilon_{1}\}|$. $S_{\max}$, $\varepsilon_{1}$ and $\varepsilon_{2}$ are some parameters which are pre-given.

With the approximate solving of APG, we obtain an approximate solution of $(\ref{equation4.1})$. Then, we obtain
\begin{eqnarray} \label{equation4.24}
\hat{z}=\left\{
\begin{array}{lcl}
l_j      &      & y_{j}^{\zeta} \leq l_j+\eta \Arrowvert y^{\zeta} \Arrowvert;\\
y_j^{\zeta}   &      & l_{j}+\eta \Arrowvert y^{\zeta} \Arrowvert <y_{j}^{\zeta} < u_{j}-\eta \Arrowvert y^{\zeta} \Arrowvert;\\
u_j       &      & y_{j}^{\zeta} \geq u_{j}- \eta \Arrowvert y^{\zeta} \Arrowvert.
\end{array} \right.
\end{eqnarray}
by filtration with $\eta$, where $\eta>0$ is a small number. The active-set of $\hat{z}$ is a good prediction to that of $\bar{z}$, which implies the steps of the following homotopy algorithm is small.

Let
\begin{eqnarray} \label{equation4.25}
 w=\left\{
\begin{array}{lcl}
\xi_{1}     &      & {\hat{z}_{j} =  l_{j}};\\
-H_{j}^T\hat{z}-f_{j}~~~          &      & {l_{j}<\hat{z}_{j}<u_{j}};\\
\xi_{2}     &      & {\hat{z}_{j} =  u_{j}},
\end{array} \right.
\end{eqnarray}
where $\xi_{1}=-\min_{j}\{H_{j}^T\hat{z}+f_{j}|\hat{z}_{j} =  l_j\}+\delta$, $\xi_{2}=-\max_{j}\{H_{j}^T\hat{z}+f_{j}|\hat{z}_{j} =  u_j\}-\delta$ and $\delta>0$.
So $\hat{z}$ is the solution of
\begin{equation} \label{equation4.26}
\begin{array}{l}
\min ~~ \frac{1}{2}z^{T}Hz+(f+w)^{T}z \\	
\rm{s.~t.}~~~\emph{l}\leq\emph{z}\leq \emph{u}.
\end{array}
\end{equation}

The linear homotopy between the objective function of $(\ref{equation4.1})$ and $(\ref{equation4.26})$ is
\begin{equation} \label{equation4.27}
\begin{array}{l}
h(t,z)= \frac{1}{2}z^{T}Hz+(f+tw)^{T}z~,~t\in[0,1].
\end{array}
\end{equation}
So we can obtain the solution of $(\ref{equation4.1})$ by solving the PQP problem
\begin{equation}\label{equation4.28}
\begin{array}{l}
\min ~~ h(t,z)=\frac{1}{2}z^{T}Hz+(f+tw)^{T}z \\	
\rm{s.~t.}~~~\emph{l}\leq \emph{z}\leq \emph{u}.
\end{array}
\end{equation}

\textbf{$\bullet$ Homotopy tracking.} $(\ref{equation4.28})$ is a special case of (\ref{pqp}), so we use PAS algorithm to track the solution path of $(\ref{equation4.28})$. Moreover, since $(\ref{equation4.28})$ is simple on constraints, we do not need to iterate the multipliers in the tracking steps, that, we simplify the PAS method to solve (\ref{equation4.28}) as follows.

Let $z(t)$, $t\in [0,1]$ be a vector-function of $t$ denoting the solution path of $(\ref{equation4.28})$ and suppose $z(t)$ is linear in $M$ intervals respectively, moreover, set $t_{0}=1, t_{M}=0$. Let $(t_{i},t_{i-1})$, $i=1,...,M$ denote the intervals, in which, $z(t)$ is linear. Let $\mathcal{A}_{m,i}=\{j|H_j z(t)+f_{j}+tw_j=0|t\in(t_{i},t_{i-1})\}$ denote the working set, $\mathcal{A}_{l,i}=\{j|H_j z(t)+f_{j}+tw_j>0|t\in(t_{i},t_{i-1})\}$ denote the lower active-set and $\mathcal{A}_{u,i}=\{j|H_j z(t)+f_{j}+tw_j<0|t\in(t_{i},t_{i-1})\}$ denote the upper active set accordingly.

\noindent
{\bf Proposition 1}\label{th:unirep}
 For any $i\in\{1,...,M\}$, there exists unique triple $(\mathcal{S}_{l,i}, \mathcal{S}_{m,i}, \mathcal{S}_{u,i})$  such that  $\mathcal{S}_{l,i}\cup \mathcal{S}_{m,i}\cup \mathcal{S}_{u,i}=\{1,...,n\}$ and
\begin{eqnarray}\label{equation4.8a}
&&z_{\mathcal{S}_{m,i}}(t) =-H_{\mathcal{S}_{m,i}\mathcal{S}_{m,i}}^{-1}(H_{\mathcal{S}_{m,i}\mathcal{S}_{u,i}}u_{\mathcal{S}_{u,i}}+H_{\mathcal{S}_{m,i}\mathcal{S}_{l,i}}l_{\mathcal{S}_{l,i}}+f_{\mathcal{S}_{m,i}}+tw_{\mathcal{S}_{m,i}}),\hspace{0.5cm}\\	 \label{equation4.8a2}
&&l_{\mathcal{S}_{m,i}}\leq z_{\mathcal{S}_{m,i}}(t)\leq u_{\mathcal{S}_{m,i}},\\\label{equation4.8b}
&&z_{\mathcal{S}_{l,i}}(t)=l_{\mathcal{S}_{l,i}},\\ \label{equation4.8c}
&&H^{T}_{\mathcal{S}_{l,i}}z(t)+f_{\mathcal{S}_{l,i}}+tw_{\mathcal{S}_{l,i}}>0, \\ \label{equation4.8d}
&&z_{\mathcal{S}_{u,i}}(t)=u_{\mathcal{S}_{u,i}},\\ \label{equation4.8e}
&&H^{T}_{\mathcal{S}_{u,i}}z(t)+f_{\mathcal{S}_{u,i}}+tw_{\mathcal{S}_{u,i}}<0,
\end{eqnarray}
for any $t\in(t_{i},t_{i-1})$, where  $H_{\mathcal{S}_{m,i}\mathcal{S}_{m,i}}$ denotes the sub-matrices of $H$ with appropriate rows and columns.

It is obvious that $(\mathcal{A}_{l,i}, \mathcal{A}_{m,i}, \mathcal{A}_{u,i})$ satisfies $(\ref{equation4.8a})$-$(\ref{equation4.8e})$ by the optimality conditions $(\ref{equation2.2a})$-$(\ref{equation2.2c})$.

Now, assume that there exists another classification $(\mathcal{B}_{l,i}, \mathcal{B}_{m,i}, \mathcal{B}_{u,i})$ that satisfies $(\ref{equation4.8a})$-$(\ref{equation4.8e})$. Since (\ref{equation4.28}) is strictly convex, for any $t$, the solution of  (\ref{equation4.28}) is unique. we have $\mathcal{A}_{l,i}=\mathcal{B}_{l,i}$ from (\ref{equation4.8b}) and (\ref{equation4.8c}), and $\mathcal{A}_{u,i}=\mathcal{B}_{u,i}$ from (\ref{equation4.8d}) and (\ref{equation4.8e}). Hence, $\mathcal{A}_{m,i}=\mathcal{B}_{m,i}$.

We start the homotopy algorithm with $t_0=1$, $i=1$, $z(t_{0})=\hat{z}$, $\mathcal{A}_{l,1}=\mathcal{A}_{l}(\hat{z})$, $\mathcal{A}_{m,1}=\mathcal{A}_{m}(\hat{z})$ and $\mathcal{A}_{u,1}=\mathcal{A}_{u}(\hat{z})$. By induction, we need to calculate $t_{i} $ and update the classification $(\mathcal{A}_{l,i+1}, \mathcal{A}_{m,i+1}, \mathcal{A}_{u,i+1})$ in the $(i+1)$-th interval, $i=1,2,...,M$.

According to the optimality conditions $(\ref{equation2.2a})$-$(\ref{equation2.2c})$, $z(t)$ satisfies (\ref{equation4.8a})-(\ref{equation4.8e}) in the $i$-th interval. Continue to decrease $t$ from $t_{i-1}$ until one of the following events occurs.
\begin{enumerate}[\rm{(}i)]
\item There exists $j\in \mathcal{A}_{m,i}$ such that $z_{j}(t)=l_{j}$.
\item There exists $j\in \mathcal{A}_{m,i}$ such that $z_{j}(t)=u_{j}$.
\item There exists $j\in \mathcal{A}_{l,i}$ such that $H_{j\mathcal{A}_{m,i}}z_{\mathcal{A}_{m,i}}(t)+H_{j\mathcal{A}_{l,i}}l_{\mathcal{A}_{l,i}}+H_{j\mathcal{A}_{u,i}}u_{\mathcal{A}_{u,i}}+(f_{j}+tw_{j})=0$.
\item There exists $j\in \mathcal{A}_{u,i}$ such that $H_{j\mathcal{A}_{m,i}}z_{\mathcal{A}_{m,i}}(t)+H_{j\mathcal{A}_{l,i}}l_{\mathcal{A}_{l,i}}+H_{j\mathcal{A}_{u,i}}u_{\mathcal{A}_{u,i}}+(f_{j}+tw_{j})=0$.
\end{enumerate}
Furthermore, according to $($i$)$-$($iv$)$, we define
\begin{eqnarray}\label{equation4.10a}
\hat{j}_{l}^{i}&=& \arg \max_{j} \{\frac{\mu^{i}_{j}-l_{j}}{\nu^{i}_{j}}< t_{i-1}|j\in \mathcal{A}_{m,i} ~and~ \nu^{i}_{j}<0\},\\  \label{equation4.10b}
\hat{j}_{u}^{i}&=&\arg\max_{j} \{\frac{\mu^{i}_{j}-u_{j}} {\nu^{i}_{j}}< t_{i-1}|j\in \mathcal{A}_{m,i} ~and~ \nu^{i}_{j}>0\},\\ \label{equation4.10c}
\tilde{j}_{l}^{i}&=& \arg \max_{j}  \{\frac{\vartheta^{i}_{j}}{\theta^{i}_{j}}< t_{i-1}|j\in \mathcal{A}_{l,i}  ~and~ \theta^{i}_{j}<0\},\\  \label{equation4.10d}
\tilde{j}_{u}^{i}&=&\arg\max_{j}  \{\frac{\varphi^{i}_{j}}{\phi^{i}_{j}}< t_{i-1}|j\in \mathcal{A}_{u,i}  ~and~ \phi^{i}_{j}>0\},
\end{eqnarray}
where
\begin{eqnarray}\label{equation4.11a}
\mu^{i}&=& -H_{\mathcal{A}_{m,i}\mathcal{A}_{m,i}}^{-1} (H_{\mathcal{A}_{m,i}\mathcal{A}_{l,i}}l_{\mathcal{A}_{l,i}}+H_{\mathcal{A}_{m,i}\mathcal{A}_{u,i}}u_{\mathcal{A}_{u,i}}+f_{\mathcal{A}_{m,i}}),\\  \label{equation4.11b}
\nu^{i}&=& H_{\mathcal{A}_{m,i}\mathcal{A}_{m,i}}^{-1} w_{\mathcal{A}_{m,i}},\\ \label{equation4.11c}
\vartheta^{i}&=&H_{\mathcal{A}_{l,i}\mathcal{A}_{m,i}}\mu^{i}+H_{\mathcal{A}_{l,i}\mathcal{A}_{l,i}}l_{\mathcal{A}_{l,i}}+H_{\mathcal{A}_{l,i}\mathcal{A}_{u,i}}u_{\mathcal{A}_{u,i}}+f_{\mathcal{A}_{l,i}},\\  \label{equation4.11d}
\theta^{i}&=& H_{\mathcal{A}_{l,i}\mathcal{A}_{m,i}}\nu^{i}-w_{\mathcal{A}_{l,i}},\\ \label{equation4.11e}
\varphi^{i}&=&H_{\mathcal{A}_{u,i}\mathcal{A}_{m,i}}\mu^{i}+H_{\mathcal{A}_{u,i}\mathcal{A}_{l,i}}l_{\mathcal{A}_{l,i}}+H_{\mathcal{A}_{u,i}\mathcal{A}_{u,i}}u_{\mathcal{A}_{u,i}}+f_{\mathcal{A}_{u,i}}, \\ \label{equation4.11f}
\phi^{i}&=& H_{\mathcal{A}_{u,i}\mathcal{A}_{m,i}}\nu^{i}-w_{\mathcal{A}_{u,i}}.
\end{eqnarray}

If there exists no $\hat{j}_{l}^{i}$ satisfies (\ref{equation4.10a}), then let $\frac{\mu^{i}_{\hat{j}_{l}^{i}}-l_{\hat{j}_{l}^{i}}}{\nu^{i}_{\hat{j}_{l}^{i}}}=-\infty$, and similarly holds for $\hat{j}_{u}^{i}$, $\tilde{j}_{l}^{i}$ and $\tilde{j}_{u}^{i}$. No matter which event occurs, the classification needs to be updated, we discuss the updating strategy in five cases.

\noindent  $\textbf{Case 1:}$  $\frac{\mu^{i}_{\hat{j}_{l}^{i}}-l_{\hat{j}_{l}^{i}}}{\nu^{i}_{\hat{j}_{l}^{i}}}>\max (\frac{\mu^{i}_{\hat{j}_{u}^{i}}-u_{\hat{j}_{u}^{i}}}{\nu^{i}_{\hat{j}_{u}^{i}}}, \frac{\vartheta^{i}_{\tilde{j}_{l}^{i}}}{\theta^{i}_{\tilde{j}_{l}^{i}}}, \frac{\varphi^{i}_{\tilde{j}_{u}^{i}}}{\phi^{i}_{\tilde{j}_{u}^{i}}})$ and $\frac{\mu^{i}_{\hat{j}_{l}^{i}}-l_{\hat{j}_{l}^{i}}}{\nu^{i}_{\hat{j}_{l}^{i}}}>0$.

\noindent It means $($i$)$ occurs first. We obtain $t_{i}=\frac{\mu^{i}_{\hat{j}_{l}^{i}}-l_{\hat{j}_{l}^{i}}}{\nu^{i}_{\hat{j}_{l}^{i}}}$, $\mathcal{A}_{m,i+1}=\mathcal{A}_{m,i}\backslash \hat{j}_{l}^{i}$, $\mathcal{A}_{l,i+1}=\mathcal{A}_{l,i}\cup \hat{j}_{l}^{i}$ and $\mathcal{A}_{u,i+1}=\mathcal{A}_{u,i}$.

\noindent   $\textbf{Case 2:}$ $\frac{\mu^{i}_{\hat{j}_{u}^{i}}-u_{\hat{j}_{u}^{i}}}{\nu^{i}_{\hat{j}_{u}^{i}}}>\max (\frac{\mu^{i}_{\hat{j}_{l}^{i}}-l_{\hat{j}_{l}^{i}}}{\nu^{i}_{\hat{j}_{l}^{i}}}, \frac{\vartheta^{i}_{\tilde{j}_{l}^{i}}}{\theta^{i}_{\tilde{j}_{l}^{i}}}, \frac{\varphi^{i}_{\tilde{j}_{u}^{i}}}{\phi^{i}_{\tilde{j}_{u}^{i}}})$ and $\frac{\mu^{i}_{\hat{j}_{u}^{i}}-u_{\hat{j}_{u}^{i}}}{\nu^{i}_{\hat{j}_{u}^{i}}}>0$.

\noindent It means $($ii$)$ occurs first. We obtain $t_{i}=\frac{\mu^{i}_{\hat{j}_{u}^{i}}-u_{\hat{j}_{u}^{i}}}{\nu^{i}_{\hat{j}_{u}^{i}}}$, $\mathcal{A}_{m,i+1}=\mathcal{A}_{m,i}\backslash \hat{j}_{u}^{i}$, $\mathcal{A}_{l,i+1}=\mathcal{A}_{l,i}$ and $\mathcal{A}_{u,i+1}=\mathcal{A}_{u,i}\cup \hat{j}_{u}^{i}$.

\noindent   $\textbf{Case 3:}$ $\frac{\vartheta^{i}_{\tilde{j}_{l}^{i}}}{\theta^{i}_{\tilde{j}_{l}^{i}}}>\max (\frac{\mu^{i}_{\hat{j}_{l}^{i}}-l_{\hat{j}_{l}^{i}}}{\nu^{i}_{\hat{j}_{l}^{i}}}, \frac{\mu^{i}_{\hat{j}_{u}^{i}}-u_{\hat{j}_{u}^{i}}}{\nu^{i}_{\hat{j}_{u}^{i}}}, \frac{\varphi^{i}_{\tilde{j}_{u}^{i}}}{\phi^{i}_{\tilde{j}_{u}^{i}}})$ and $\frac{\vartheta^{i}_{\tilde{j}_{l}^{i}}}{\theta^{k,i}_{\tilde{j}_{l}^{i}}}>0$.

\noindent It means $($iii$)$ occurs first, then we obtain $t_{i}=\frac{\vartheta^{i}_{\tilde{j}_{l}^{i}}}{\theta^{i}_{\tilde{j}_{l}^{i}}}$, $\mathcal{A}_{m,i+1}=\mathcal{A}_{m,i}\cup \tilde{j}_{l}^{i}$, $\mathcal{A}_{l,i+1}=\mathcal{A}_{l,i}\backslash \tilde{j}_{l}^{i}$ and $\mathcal{A}_{u,i+1}=\mathcal{A}_{u,i}$.

\noindent   $\textbf{Case 4:}$ $\frac{\varphi^{i}_{\tilde{j}_{u}^{i}}}{\phi^{i}_{\tilde{j}_{u}^{i}}}>\max (\frac{\mu^{i}_{\hat{j}_{l}^{i}}-l_{\hat{j}_{l}^{i}}}{\nu^{i}_{\hat{j}_{l}^{i}}},\frac{\mu^{i}_{\hat{j}_{u}^{i}}-u_{\hat{j}_{u}^{i}}}{\nu^{i}_{\hat{j}_{u}^{i}}}, \frac{\vartheta^{i}_{\tilde{j}_{l}^{i}}}{\theta^{i}_{\tilde{j}_{l}^{i}}})$ and $ \frac{\varphi^{i}_{\tilde{j}_{u}^{i}}}{\phi^{i}_{\tilde{j}_{u}^{i}}}>0$.

\noindent It means $($iv$)$ occurs first, then we obtain $t_{i}=\frac{\varphi^{i}_{\tilde{j}_{u}^{i}}}{\phi^{i}_{\tilde{j}_{u}^{i}}}$, $\mathcal{A}_{m,i+1}=\mathcal{A}_{m,i}\cup \tilde{j}_{u}^{i}$, $\mathcal{A}_{l,i+1}=\mathcal{A}_{l,i}$ and $\mathcal{A}_{u,i+1}=\mathcal{A}_{u,i}\backslash \tilde{j}_{u}^{i}$.

\noindent  $\textbf{Case 5:}$  $\max (\frac{\mu^{i}_{\hat{j}_{l}^{i}}-l_{\hat{j}_{l}^{i}}}{\nu^{i}_{\hat{j}_{l}^{i}}},\frac{\mu^{i}_{\hat{j}_{u}^{i}}-u_{\hat{j}_{u}^{i}}}{\nu^{i}_{\hat{j}_{u}^{i}}}, \frac{\vartheta^{i}_{\tilde{j}_{l}^{i}}}{\theta^{i}_{\tilde{j}_{l}^{i}}}, \frac{\varphi^{i}_{\tilde{j}_{u}^{i}}}{\phi^{i}_{\tilde{j}_{u}^{i}}})\leq 0$.

\noindent In this case, the algorithm ends and we obtain
\begin{eqnarray}\label{equation4.13a}
&&z_{\mathcal{A}_{m,i}}(0) =-H_{\mathcal{A}_{m,i}\mathcal{A}_{m,i}}^{-1} (H_{\mathcal{A}_{m,i}\mathcal{A}_{u,i}}u_{\mathcal{A}_{u,i}}+H_{\mathcal{A}_{m,i}\mathcal{A}_{l,i}}l_{\mathcal{A}_{l,i}}+f_{\mathcal{A}_{m,i}}),\\	\label{equation4.13b}
&&z_{\mathcal{A}_{l,i}}(0)=l_{\mathcal{A}_{l,i}}, \\ \label{equation4.13c}
&&z_{\mathcal{A}_{u,i}}(0)=u_{\mathcal{A}_{u,i}}.
\end{eqnarray}

\textbf{$\bullet$ $\varepsilon$-precision verification and correction.} From the homotopy tracking steps, we have
\begin{eqnarray} \nonumber
&&z_{\mathcal{A}_{m,i}}(t_i)=-H_{\mathcal{A}_{m,i}\mathcal{A}_{m,i}}^{-1} (H_{\mathcal{A}_{m,i}\mathcal{A}_{u,i}}u_{\mathcal{A}_{u,i}}+H_{\mathcal{A}_{m,i}\mathcal{A}_{l,i}}l_{\mathcal{A}_{l,i}}+f_{\mathcal{A}_{m,i}}+t_{i}w_{\mathcal{A}_{m,i}}), \\ \nonumber
&&z_{\mathcal{A}_{l,i}}(t_i)=l_{\mathcal{A}_{l,i}}, \\\nonumber
&&z_{\mathcal{A}_{u,i}}(t_i)=u_{\mathcal{A}_{u,i}}.
\end{eqnarray}
However, due to the errors from the solving of the linear systems, the solution path obtained from the tracking steps above may not satisfy the optimality conditions, so we need to verify that  $z(t_i)$ satisfies  the optimality conditions:

\begin{eqnarray} \label{kkt at ti}
&&l_{\mathcal{A}_{m,i}}\leq z_{\mathcal{A}_{m,i}}(t_i)\leq u_{\mathcal{A}_{m,i}}\\\label{kkt at ti1}
&&H^T_{\mathcal{A}_{l,i}}z(t_i)+f_{\mathcal{A}_{l,i}}+t_{i}w_{\mathcal{A}_{l,i}}\geq 0\\\label{kkt at ti2}
&&H^T_{\mathcal{A}_{u,i}}z(t_i)+f_{\mathcal{A}_{u,i}}+t_{i}w_{\mathcal{A}_{u,i}}\leq 0
\end{eqnarray}
In practice, it is not necessary and may be hard to ensure (\ref{kkt at ti})-(\ref{kkt at ti2}) hold strictly, so we relax (\ref{kkt at ti})-(\ref{kkt at ti2}) by a small $\varepsilon$ as follows:
\begin{eqnarray} \label{kkt at ti3}
&&l_{\mathcal{A}_{m,i}}(t_i)-\varepsilon\leq z_{\mathcal{A}_{m,i}}(t_i)\leq u_{\mathcal{A}_{m,i}}(t_i)+\varepsilon\\\label{kkt at ti4}
&&H^T_{\mathcal{A}_{l,i}}z(t_i)+f_{\mathcal{A}_{l,i}}+t_{i}w_{\mathcal{A}_{l,i}}\geq -\varepsilon\\\label{kkt at ti5}
&&H^T_{\mathcal{A}_{u,i}}z(t_i)+f_{\mathcal{A}_{u,i}}+t_{i}w_{\mathcal{A}_{u,i}}\leq \varepsilon
\end{eqnarray}
If (\ref{kkt at ti3})-(\ref{kkt at ti5}) hold, the homotopy goes to the next step; otherwise, we correct the $\mathcal{A}_{l,i}, \mathcal{A}_{m,i}, \mathcal{A}_{u,i}$ as follows:

\textbf{Step 1}: If there exists $J\subset \mathcal{A}_{m,i}$ such that $z_J(t_i)<l_{J}-\varepsilon$, then let
$$\bar{j}=\arg\min_{j\in J} ~\{z_j(t_i)-l_{j}+\varepsilon\}$$
and $\mathcal{A}_{l,i}=\mathcal{A}_{l,i}\cup \bar{j}, \mathcal{A}_{m,i}=\mathcal{A}_{m,i}\backslash \bar{j}$, refresh $z(t_i)$ like (\ref{kkt at ti})-(\ref{kkt at ti2}) and go to \textbf{Step 1}; otherwise go to \textbf{Step 2}.

\textbf{Step 2}:  If there exists $J\subset \mathcal{A}_{m,i}$ such that $z_J(t_i)>u_{J}+\varepsilon$, then let
$$\bar{j}=\arg\max_{j\in J}~ \{z_j(t_i)-u_{j}-\varepsilon\}$$
and $\mathcal{A}_{u,i}=\mathcal{A}_{u,i}\cup \bar{j}, \mathcal{A}_{m,i}=\mathcal{A}_{m,i}\backslash \bar{j}$, refresh $z(t_i)$ like (\ref{kkt at ti})-(\ref{kkt at ti2}) and go to  \textbf{Step 1}; otherwise go to \textbf{Step 3}.

\textbf{Step 3}: If there exists $J\subset \mathcal{A}_{l,i}$ such that $H^T_{J}z(t_i)+f_{J}+t_{i}w_{J}\leq -\varepsilon$, then let
$$\bar{j}=\arg\min_{j\in J} ~\{H^T_{j}z(t_i)+f_{j}+t_{i}w_{j} +\varepsilon\}$$
and $\mathcal{A}_{l,i}=\mathcal{A}_{l,i}\backslash \bar{j}, \mathcal{A}_{m,i}=\mathcal{A}_{m,i}\cup \bar{j}$, refresh $z(t_i)$ like (\ref{kkt at ti})-(\ref{kkt at ti2}) and go to \textbf{Step 1}; otherwise go to \textbf{Step 4}.

\textbf{Step 4}: If there exists $J\subset \mathcal{A}_{u,i}$ such that $H^T_{J}z(t_i)+f_{J}+t_{i}w_{J}\geq \varepsilon$, then let
$$\bar{j}=\arg\max_{j\in J} ~\{H^T_{j}z(t_i)+f_{j}+t_{i}w_{j} -\varepsilon\}$$
and $\mathcal{A}_{u,i}=\mathcal{A}_{u,i}\backslash \bar{j}, \mathcal{A}_{m,i}=\mathcal{A}_{m,i}\cup \bar{j}$, refresh $z(t_i)$ like (\ref{kkt at ti})-(\ref{kkt at ti2}) and go to \textbf{Step 1}; otherwise end the correction.

The correction steps ensure the solution $z(t)$ satisfies the optimality conditions with $\varepsilon$-precision and can guarantee the stability of the homotopy tracking algorithm.

 At each step of the homotopy algorithm, we need to solve two symmetric positive definite linear systems of equations
\begin{eqnarray} \label{equation4.14a}
H_{\mathcal{A}_{m,i}\mathcal{A}_{m,i}}\mu^{i}&=&-H_{\mathcal{A}_{m,i}\mathcal{A}_{l,i}}l_{\mathcal{A}_{l,i}}-H_{\mathcal{A}_{m,i}\mathcal{A}_{u,i}}u_{\mathcal{A}_{u,i}}-f_{\mathcal{A}_{m,i}},\\  \label{equation4.14b}
H_{\mathcal{A}_{m,i}\mathcal{A}_{m,i}}\nu^{i}&=&  w_{\mathcal{A}_{m,i}},
\end{eqnarray}
and do matrix-vector multiplications in $(\ref{equation4.11c})$-$(\ref{equation4.11f})$. From (\ref{equation4.14a}) and  (\ref{equation4.14b}), we see the homotopy algorithm takes small computation at each step if $|\mathcal{A}_{m,i}|$ is small, that is, the number of the free variables is small, which implies the homotopy algorithm would be efficient for solving BQPs in SVM.

Since $H_{\mathcal{A}_{m,i}\mathcal{A}_{m,i}}$ is positive definite, we apply Cholesky decomposition method for (\ref{equation4.14a}) and (\ref{equation4.14b}). Moreover, since $\mathcal{A}_{m,i}$ changes only one index at each step, we implement the   Cholesky decomposition technique   \cite{ferreau2006online,ferreau2014qpoases} for (\ref{equation4.14a}) and  (\ref{equation4.14b}).

We conclude the framework of the Homotopy algorithm as Algorithm \ref{algorithm1}.
\vspace{-0.5cm}

\begin{algorithm}[!h]
    \caption{Homtopy algorithm for strictly convex BQP problem}

    \begin{algorithmic}[]
    \REQUIRE~$z^{0}$;
    \ENSURE~$\bar{z}$.

       \STATE Approximately solve $(\ref{equation4.1})$ with  APG algorithm until $(\ref{equation4.21a})$ or $(\ref{equation4.21b})$ is satisfied;

       \STATE Obtain a warm start $\hat{z}$ like $(\ref{equation4.24})$ and generate a PQP problem $(\ref{equation4.28})$;

       \STATE Track the solution path of $(\ref{equation4.28})$ from $t=1$ to $t=0$ and obtain $\bar{z}$ in a form like (\ref{equation4.13a})-(\ref{equation4.13c}).

    \end{algorithmic} \label{algorithm1}

\end{algorithm}

\subsection{Practical proximal homotopy algorithms for non-convex BQP}
Note that it is not necessary to compute the exact solution of $(\ref{equation3.1})$ in the previous iterations of PP, so we directly go to the next iteration when $y^{l}$ makes $q(x)$ decrease enough, that is, $y^{l}$ satisfies
\begin{equation} \label{suffidient_decrease}
\begin{array}{l}
q(y^{l})\leq q(x^{k})-f_{\varepsilon},
\end{array}
\end{equation}
where $f_{\varepsilon}>0$ is given. Then we have practical PP-Hom and APP-Hom algorithms as follows
%%%%%%%%%%%%%%%%%%%%%%%%%%%%%%%%%%%%%%%%%%%%%%%%%%%%%%%%%%%%%%%%%%%%%%%%%%%%%%%%%%

\begin{algorithm}[!h]
    \caption{PP-Hom algorithm for non-convex BQP}

    \begin{algorithmic}[]
    \REQUIRE ~~                         %Parameters Input
    $x^{0}$, $tol$, $f_{\varepsilon}$;
    \ENSURE ~~                         %Output
    $x^{k+1}$;
    \WHILE{$\|x^k-x^{k+1}\|>tol$}

       \STATE Approximately solve $(\ref{equation3.1})$ with APG until $(\ref{equation4.21a})$ or $(\ref{equation4.21b})$ is satisfied;

       \IF{$q(y^{l})\leq q(x^{k})-f_{\varepsilon}$}
       \STATE $x^{k+1}=y^{l}$;
        \ELSE
       \STATE Obtain a warm start $\hat{z}$ like $(\ref{equation4.24})$ and generate a PQP problem $(\ref{equation4.28})$;
        \STATE Track the solution path of $(\ref{equation4.28})$ from $t=1$ to $t=0$ and return $x^{k+1}$ in a form like (\ref{equation4.13a})-(\ref{equation4.13c}).
        \ENDIF

    \ENDWHILE
    \end{algorithmic} \label{algorithm2}

\end{algorithm}

\vspace{-1cm}

%%%%%%%%%%%%%%%%%%%%%%%%%%%%%%%%%%%%%%%%%%%%%%%%%%%%%%%%%%%%%%%%%%%%%%%%%%%%%%%%
\begin{algorithm}[!h]
    \caption{APP-Hom algorithm for non-convex BQP}

    \begin{algorithmic}[]
    \REQUIRE ~~                         %Parameters Input
    $x^{0}$, $tol$, $f_{\varepsilon}$;
    \ENSURE ~~                          %Output
    $x^{k+1}$;
    \WHILE{$\|x^k-x^{k+1}\|>tol$}
     \IF{$\omega_k$ satisfies (\ref{switch condition})}

       \STATE Approximately solve $(\ref{accelerate PPA})$ with APG  until $(\ref{equation4.21a})$ or $(\ref{equation4.21b})$ is satisfied;
       \ELSE
        \STATE Approximately solve $(\ref{equation3.1})$ with APG until $(\ref{equation4.21a})$ or $(\ref{equation4.21b})$ is satisfied;
     \ENDIF
       \IF{$q(y^{l})\leq q(x^{k})-f_{\varepsilon}$}
       \STATE $x^{k+1}=y^{l}$;
        \ELSE

       \STATE Obtain a warm start $\hat{z}$ like $(\ref{equation4.24})$ and generate a PQP problem $(\ref{equation4.28})$;
        \STATE Track the solution path of $(\ref{equation4.28})$ from $t=1$ to $t=0$ and return $x^{k+1}$ in a form like (\ref{equation4.13a})-(\ref{equation4.13c}).
        \ENDIF

    \ENDWHILE
    \end{algorithmic} \label{algorithm3}

\end{algorithm}

\section{Numerical experiments}
In this section, we present numerical results obtained from the implementations of our algorithms described above. The numerical experiments were performed by Matlab 8.1 programming platform (R2013a) running on a machine with Windows 8 Operation System, Intel(R) Core(TM)i7 CPU 4790 3.60GHz processor and 32GB of memory. For convenience, we use TRR to denote the reflective Newton method (the ``trust-region-reflective'' algorithm in Matlab $quadprog$ solver).
In our experiments, we terminated PP and APP when $\|x^k-x^{k-1}\|<10^{-11}$. Moreover, let
\begin{eqnarray} \label{gradx}
g(x)=\left\{
\begin{array}{lcl}
Q_{i}^{T}x+r_i,      &      & i\in\mathcal{A}_{m}(x);\\
0,  &      & else.
\end{array} \right.
\end{eqnarray}

\subsection{Nonnegative least-squares problems}

Nonnegative least-squares problem \cite{lawson1995solving}
 \begin{equation} \label{nnls}
\begin{array}{l}
\min_{x\in \mathbb{R}^{n}}~\frac{1}{2}\Arrowvert Ax-b\Arrowvert^{2}~~~~s.t.~~ x\geq 0,
\end{array}
\end{equation}
is a classical problem in scientific computing. Applications include image restoration \cite{nagy2000enforcing}, non-negative matrix factorization \cite{berry2007algorithms}, etc. Many algorithms can be applied to solve this kind of problems, e.g., LBFGS-B, TRR and FNNLS, where FNNLS \cite{bro1997fast} is Bro and Jong's improved implementation of the Lawson-Hanson NNLS procedure \cite{lawson1995solving}.

$\bullet$ \textbf{Random NNLS problems:} Given $m>n$, we randomly generated NNLS problems with Matlab codes:

\qquad $A$={\rm sprandn}($m$,$n$,$\rho$); $\bar{x}$={\rm randn}($n$,1); $\bar{x}=\bar{x}.*(\bar{x}>0)$; $b=A*\bar{x}$;

The data was divided into two cases: dense NNLS problems (NNLS-D) and sparse NNLS problems (NNLS-S).

We first compared the homotopy algorithm with LBFGS-B, ASA, MINQ8 \cite{huyer2018minq8} and FNNLS on solving NNLS-D and NNLS-S. We use the residual $\|Ax-b\|$ and $\|g(x)\|$ to measure the precision.
The results are shown in Table \ref{table3} and \ref{table4}. It shows that the homotopy algorithm outperforms LBFGS-B, ASA, MINQ8 and FNNLS. FNNLS is not suitable for large-scale problems.
\begin{table}[!htb]
 \caption{\scriptsize{Solving random dense NNLS problems: NNLS-D }}

    \label{table3}
     \tiny

    \begin{tabular}{p{1.9cm}p{1.5cm}p{1.5cm}p{1.5cm}p{1.5cm}p{1.5cm}}
    \toprule
    Methods& NNLS-D1 &NNLS-D2 &NNLS-D3&NNLS-D4&NNLS-D5\\
    \midrule
   $m\times n$&1000$\times$800  &2000$\times$500 &5000$\times$4000&8000$\times$7000&15000$\times$6000\\
    \midrule
   Homtopy&0.11s  &0.05s &1.10s&4.34s&1.67s\\
   $\| g(x)\|$&1.69E-11 &1.86E-11 &4.19E-10&1.24E-09&2.98E-09\\
   $\|Ax-b\|$&7.26E-13  &5.65E-13 &6.81E-12&1.56E-11&2.14E-11\\
   \midrule
   LBFGS-B&0.09s  &0.06s &1.18s&4.07s&3.19s\\
   $\| g(x)\|$&1.82E-04 &2.38E-04 &3.14E-03&4.25E-03&1.33E-03\\
   $\|Ax-b\|$&1.46E-05  &7.69E-06&1.02E-04&1.16E-04&2.01E-05\\
   \midrule
   FNNLS&0.55s  &0.07s &91.17s&655.77s&330.41s\\
   $\| g(x)\|$&2.11E-11 &2.00E-11 &5.60E-10&1.49E-09&2.29E-09\\
   $\|Ax-b\|$&7.44E-13  &5.99E-13 &8.71E-12&2.00E-11&2.11E-11\\
   \midrule
   MINQ8&0.31s  &0.11s &24.75s&175.16s&133.22s\\
   $\| g(x)\|$&1.15E-05 &7.12E-05 &7.80E-04&6.59E-04&1.21E-04\\
   $\|Ax-b\|$&1.21E-06  &7.11E-13 &2.61E-05&1.94E-05&5.21E-05\\
   \midrule
    ASA&2.41s  &0.13s &9.81s&44.61s&10.52s\\
   $\| g(x)\|$&7.56E-08 &4.73E-09 &8.51E-08&8.49E-08&6.96E-08\\
   $\|Ax-b\|$&6.47E-08  &1.50E-09 &6.15E-09&3.93E-09&1.41E-09\\
    \bottomrule
    \end{tabular}
\vspace{-0.2cm}
\end{table}

\begin{table}[!htb]
 \caption{\scriptsize{Solving random sparse NNLS problems: NNLS-S }}

    \label{table4}
    \tiny

    \begin{tabular}{p{1.9cm}p{1.5cm}p{1.5cm}p{1.5cm}p{1.5cm}p{1.5cm}}
    \toprule
    Methods& NNLS-S1 &NNLS-S2 &NNLS-S3&NNLS-S4&NNLS-S5\\
     \midrule
   $m\times n$  &10000$\times$9000 &20000$\times$16000 &30000$\times$15000 &50000$\times$30000 &80000$\times$20000\\
   sparsity&0.999&0.999&0.9998&0.9998&0.9998\\
    \midrule
   Homtopy&1.97s  &4.82s &1.41s&18.24s&3.89s\\
   $\| g(x)\|$&1.45E-12 &5.22E-12 &7.65E-12&3.84E-12&4.39E-12\\
   $\|Ax-b\|$&5.13E-13  &1.43E-12 &1.51E-12&1.28E-12&1.07E-12\\
   \midrule
   LBFGS-B&4.54s  &4.43s &2.54s&41.18s&5.50s\\
   $\| g(x)\|$&1.89E-05 &2.27E-05 &8.16E-06&2.27E-05&1.75E-05\\
   $\|Ax-b\|$&1.10E-04  &1.04E-05&3.25E-06&5.02E-05&8.91E-06\\
   \midrule
   FNNLS&897.12s  &6902.58s &5201.17s&Hours&Hours\\
   $\| g(x)\|$&1.30E-12 &2.00E-11 &5.54E-12&-&-\\
   $\|Ax-b\|$&5.01E-13  &5.99E-13 &3.17E-12&-&-\\
   \midrule
   MINQ8&588.31s  &3122.41s &2987.35s&Hours&Hours\\
   $\| g(x)\|$&7.00E-03 &6.32E-03 &1.77E-03&5.51E-03&4.00E-03\\
   $\|Ax-b\|$&5.21E-03  &3.44E-03 &3.98E-03&2.88E-03&3.11E-03\\
   \midrule
    ASA&12.81s  &13.98s &618.01s&196.72s&4.66s\\
   $\| g(x)\|$&2.21E-08 &1.08E-08 &1.30E-08&8.49E-08&9.74E-09\\
   $\|Ax-b\|$&1.07E-07  &4.86E-06 &2.15E-07&3.93E-09&1.15E-08\\
    \bottomrule
    \end{tabular}
\vspace{-0.0cm}

\end{table}

$\bullet$ \textbf{Image debluring:}
Image deblurring \cite{benvenuto2009nonnegative, hanke2000quasi} is a linear inverse problem, which discrete form is
$$Ax+\eta =y,$$
where $A\in \mathbb{R}^{N^2\times N^2}$ is a large ill-conditioned matrix representing the blurring phenomena, $\eta$ is modeling noise. The vector $x$ represents the unknown true image, $y$ is the blurred-noisy copy of $x$.

A typical model for deblurring is NNLS:
 \begin{equation} \label{deblur}
\begin{array}{l}
\min_{x\in \mathbb{R}^{n}}~\frac{1}{2}\Arrowvert Ax-y\Arrowvert^{2}~~~~s.t.~~ x\geq 0.
\end{array}
\end{equation}
Since problem (\ref{deblur}) is ill-posed, it should be regularized. We add the Tikhonov regularization $\frac{\beta}{2}\|x\|^2$ to the objective, that is,
 \begin{equation} \label{regu deblur}
\begin{array}{l}
\min_{x\in \mathbb{R}^{n}}~\frac{1}{2}\Arrowvert Ax-y\Arrowvert^{2}+\frac{\beta}{2}\|x\|^2~~~~s.t.~~ x\geq 0.
\end{array}
\end{equation}
Then for any $y$, the solution of (\ref{regu deblur}) is unique.

We obtained the satellite image from \cite{hanke2000quasi} and the star image from \cite{nagy1998restoring}. The black-white grid image was generated by ourselves. All of them are in size $256\times 256$. We used the motion blur kernel to blur the satellite image and the circular averaging blur kernel to blur the star image. The black-white grid image was blurred by Gaussian blur kernel. The noise $\eta$ was obtained from Gaussian noise with intensity 0.01, that is, $\eta=0.01*{\rm randn}(N^2,1)$. We used the homotopy, LBFGS-B, ASA, MINQ8 and TRR to solve (\ref{regu deblur}) with $\beta=10^{-4}$, and we did not present the results of FNNLS for it does not scale on these problems. The running time and precision are listed in Table \ref{table5}. The homotopy algorithm takes less time and obtains higher-precision solution.

\begin{figure}[!htb]
  \centering
  \vspace{-3.2cm}
  \includegraphics[width=11cm]{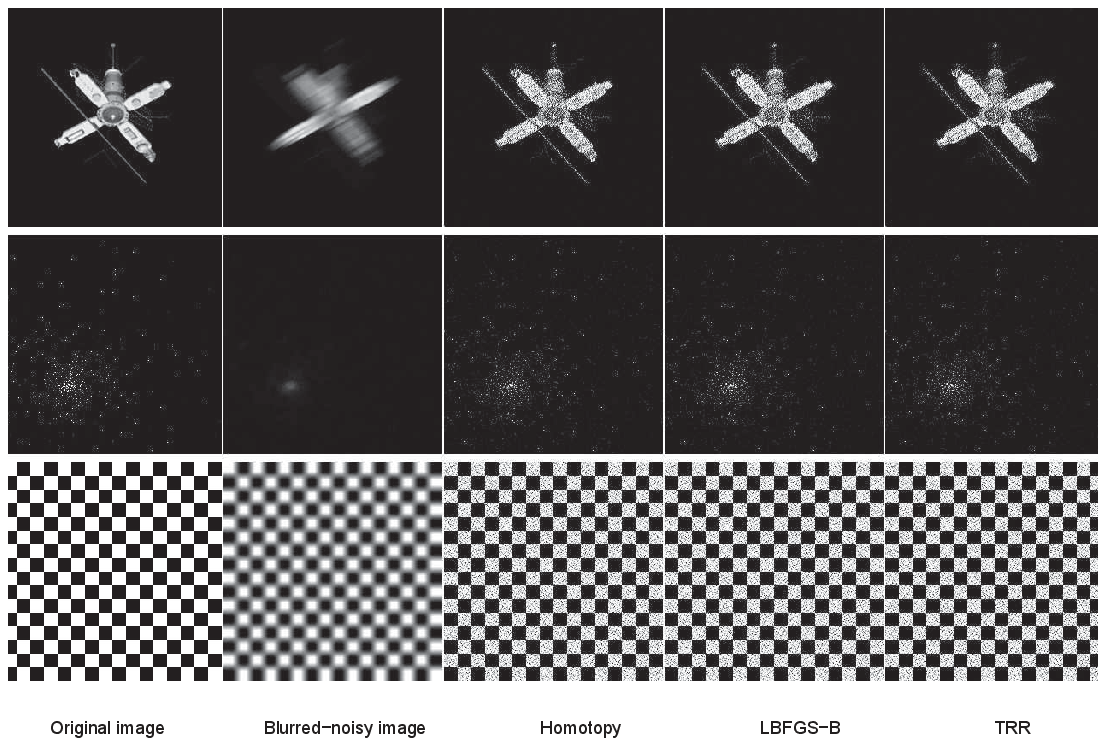}\\
  \vspace{-3.2cm}
  \caption{Deblurring using the homotopy, LBFGS-B, TRR, from top to bottom, a satellite, a star and a white-black grid image.}
\label{figure1}
\end{figure}
\vspace{-0.2cm}

\begin{table}[!htb]
 \caption{\scriptsize{Deblurring with $\beta=10^{-4}$} and solving by the homotopy, LBFGS-B, TRR.}

    \label{table5}
     \tiny

    \begin{tabular}{p{2.6cm}p{2.5cm}p{2.5cm}p{2.5cm}p{2.5cm}}
    \toprule
    Methods& Satellite image&Star image&Black-white grid image\\
   \midrule
   Homotopy&21.76s  &34.94s &51.75s\\
   $\| g(x)\|$&2.46E-14 &2.29E-14 &2.98E-14\\
   \midrule
   LBFGS-B&59.86s  &73.21s &79.90s\\
   $\| g(x)\|$&2.22E-07 &5.08E-08 &3.50E-06\\
   \midrule
   TRR&60.39s  &165.77s &90.62s\\
   $\| g(x)\|$&3.67E-04 &1.18E-04 &1.35E-03\\
   \midrule
   MINQ8&1844.76s  &2311.42s &3122.21s\\
   $\| g(x)\|$&1.32E-05 &4.21E-05 &6.12E-05\\
   \midrule
   ASA&94.50s  &95.20s &140.36s\\
   $\| g(x)\|$&1.45E-07 &1.31E-07 &3.07E-07\\
    \bottomrule
    \end{tabular}

\vspace{-0.1cm}
\end{table}

\subsection{General BQP}
In this subsection, we tested our algorithms on solving general BQPs including random non-convex BQPs, support vector machine (SVM) optimization problems and three partial differential optimization problems: the obstacle problem \cite{philippe1978finite}, the elastic-plastic torsion problem \cite{glowinski1985numerical}, and the journal bearing problem \cite{capriz1983free, cimatti1976problem}.

$\bullet$ \textbf{Non-convex BQP.} In this part, we tested PP-Hom and APP-Hom on solving non-convex BQPs: the first database was generated by ourselves; the second database was downloaded from website\footnotemark[1]\footnotetext[1]{http://www.minlp.com/nlp-and-minlp-test-problems}, which was generated by Vandenbussche and Nemhauser \cite{vandenbussche2005branch}. We generated data in dense (NCBQP-D) and sparse (NCBQP-S) cases. The Matlab codes are as follows:

~~~~$B={\rm sprandn}(n,n,\rho);Q=B'+B+\lambda*{\rm speye}(n,n)$;$r$={\rm randn}($n$,1);

~~~~$l$=zeros($n$,1);$u$=10*ones($n$,1);

\begin{table}[!htb]
 \caption{\scriptsize{Random dense non-convex BQPs: NCBQP-D with $\lambda=10$.}}

    \label{table8}
    \tiny

    \begin{tabular}{p{1.8cm}p{1.5cm}p{1.5cm}p{1.5cm}p{1.5cm}p{1.5cm}}
    \toprule
    Methods& NCBQP-D1 &NCBQP-D2 &NCBQP-D3&NCBQP-D4&NCBQP-D5\\
    \midrule
    $n$&1000  &2000 &3000&4000&5000\\
    \midrule
   PP-Hom& 2.87s&28.44s  &93.22s &127.44s&146.44s\\
   $\| g(x)\|$&7.84E-10 &1.16E-09&1.44E-09 &1.38E-09&1.38E-09\\
   $\lambda_{\min}(Q_{\bar{\mathcal{A}}_{m}\bar{\mathcal{A}}_{m}})$&0.241 &0.053&0.189 &0.134&0.245\\
   \midrule
   APP-Hom&0.41s  &2.78s &8.09s&16.55s&23.39s\\
   $\| g(x)\|$&6.38e-10 &1.06E-09 &1.25E-09&1.40E-09&1.58E-09\\
   $\lambda_{\min}(Q_{\bar{\mathcal{A}}_{m}\bar{\mathcal{A}}_{m}})$&0.495 &0.118&0.201 &0.123&0.133\\
   \midrule
   LBFGS-B&0.91s  &4.83s &16.43s&26.19s&60.25s\\
   $\| g(x)\|$&4.05E-05 &7.63E-05 &2.73E-04&1.28E-03&1.28E-03\\
   \midrule
   TRR&1.43s  &6.70s&19.60s&65.21s&89.67s\\
   $\| g(x)\|$&3.39E-06 &9.81E-06 &4.30E-06&1.58E-05&5.11E-06\\
   \midrule
   MINQ8&0.61s  &2.96s &9.95s&20.33s&28.22\\
   $\| g(x)\|$&2.80E-05 &2.68E-05 &2.55E-05&1.53E-05&3.23E-05\\

   \midrule
    ASA&0.51s  &2.88s &9.02s&17.33s&24.11s\\
   $\| g(x)\|$&2.32E-08 &1.43E-08 &2.31E-08&3.55E-08&2.86E-09\\

    \bottomrule
    \end{tabular}

\end{table}

\begin{table}[!htb]
 \caption{\scriptsize{Random sparse non-convex BQPs: NCBQP-S with $\lambda=1$. }}

    \label{table9}
    \tiny

    \begin{tabular}{p{1.8cm}p{1.5cm}p{1.5cm}p{1.5cm}p{1.5cm}p{1.5cm}}
    \toprule
    Methods& NCBQP-S1 &NCBQP-S2 &NCBQP-S3&NCBQP-S4&NCBQP-S5\\
    \midrule
   $ n$  &3000 &5000 &8000 &10000 &15000\\
   sparsity&0.99&0.99&0.99&0.99&0.99\\
    \midrule

   PP-Hom& 1.37s&3.10s  &12.81s &24.67s&67.45s\\
   $\| g(x)\|$&1.13E-10 &1.41E-10&1.82E-10 &2.01E-10&1.78E-10\\
$\lambda_{\min}(Q_{\bar{\mathcal{A}}_{m}\bar{\mathcal{A}}_{m}})$&2.65 &1.89&1.11 &1.51&0.44\\
   \midrule
   APP-Hom&0.25s &0.64s &1.68s&2.11s&6.53s\\
   $\| g(x)\|$&1.09E-10 &1.23E-10 &1.70E-10&1.44E-10&1.33E-10\\
$\lambda_{\min}(Q_{\bar{\mathcal{A}}_{m}\bar{\mathcal{A}}_{m}})$&3.42 &1.33&0.86 &1.43&0.24\\
   \midrule
   LBFGS-B&0.60s  &1.52s &5.30s&6.68s&25.56s\\
   $\| g(x)\|$&1.04E-05 &1.13E-05 &3.61E-05&2.25E-04&1.51E-06\\
   \midrule
   TRR&0.44s  &0.99s&2.03s&3.02s&9.16s\\
   $\| g(x)\|$&2.24E-06 &1.72E-06 &1.47E-06&2.00E-07&9.16E-05\\
   \midrule
   MINQ8&1.21s  &4.86s &9.95s&21.72s&28.22\\
   $\| g(x)\|$&2.80E-05 &3.38E-05 &3.97E-05&1.53E-05&3.23E-05\\
   \midrule
    ASA&0.31s  &1.02s &2.11s&3.45s&7.33s\\
   $\| g(x)\|$&2.66E-09 &2.32E-09 &6.45E-09&8.12E-09&3.66E-09\\

    \bottomrule
    \end{tabular}

\end{table}

\begin{figure}[!htb]
  \centering
  \vspace{-3cm}
  \includegraphics[width=7cm]{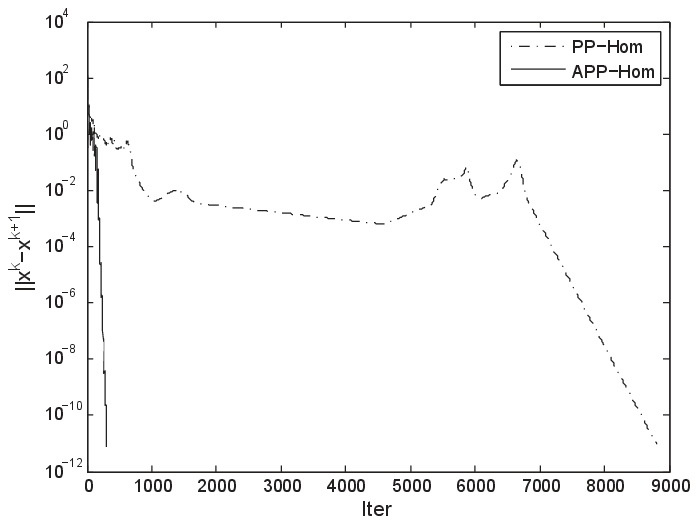}\hspace{-2.5cm}
  \includegraphics[width=7cm]{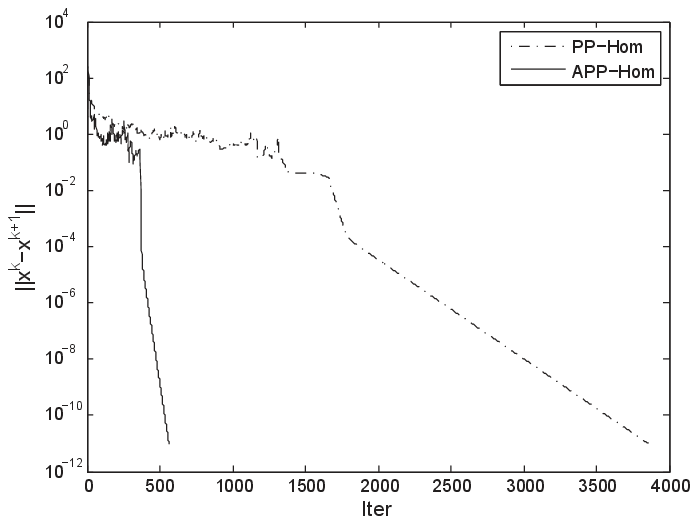}\\
  \vspace{-3cm}
  \caption{Iterative points distance for NCBQP-D3 (left) and  NCBQP-S5 (right).}\label{figure2}

\end{figure}
%%%%%%%%%%%%%%%%%%%%%%%%%%%%%%%%%%%%%%%%%%%%%%%%%%%%%%%%%%%%%%%%%%%%%%%%%%%%%%%%%%
\begin{table}[!htb]
 \caption{\scriptsize{PP-Hom and APP-Hom solving the non-convex BQPs from \cite{vandenbussche2005branch}}}

    \label{table10}
     \tiny

    \begin{tabular}{p{1.8cm}p{0.5cm}p{1cm}p{0.5cm}p{1.2cm}p{0.5cm}p{1cm}p{0.5cm}p{1.2cm}p{0.1cm}}
    \toprule
    \multirow{2}{*}{Problems}&\multirow{2}{*}{n}& \multicolumn{3}{c}{PP-Hom} && \multicolumn{3}{c}{APP-Hom} \\
    \cline{3-5}\cline{7-9}
       &&Time&Iter&$\| g(x)\|$&&Time&Iter&$\| g(x)\|$\\
    \hline
    spar020$\_100\_1$ &20&0.02s&92&2.21E-09&  &0.01s&36&6.10E-10\\
    spar020$\_100\_2$ &20&0.00s&4&0.00E+00&  &0.00s&4&0.00E+00\\
    spar020$\_100\_3$ &20&0.00s&17&1.65E-09&  &0.00s&8&3.13E-13\\
    spar030$\_60\_1$ &30&0.01s&39&0.00E+00&  &0.00s&30&0.00E+00\\
    spar030$\_60\_2$ &30&0.01s&27&0.00E+00&  &0.00s&20&0.00E+00\\
    spar030$\_60\_3$ &30&0.02s&113&0.00E+00&  &0.00s&23&0.00E+00\\
    spar030$\_70\_1$ &30&0.01s&136&1.87E-09&  &0.00s&46&2.00E-14\\
    spar030$\_70\_2$ &30&0.01s&18&0.00E+00&  &0.00s&16&0.00E+00\\
    spar030$\_70\_3$ &30&0.02s&113&1.78E-09&  &0.01s&63&1.01E-09\\
    spar030$\_80\_1$ &30&0.01s&29&2.32E-09&  &0.00s&21&0.00E+00\\
    spar030$\_80\_2$ &30&0.02s&81&0.00E+00&  &0.00s&14&0.00E+00\\
    spar030$\_80\_3$ &30&0.02s&105&2.33E-09&  &0.00s&36&4.30E-11\\
    spar030$\_90\_1$ &30&0.02s&122&2.43E-09&  &0.01s&42&1.01E-09\\
    spar030$\_90\_2$ &30&0.23s&1277&0.00E+00&  &0.01s&56&7.11E-15\\
    spar030$\_90\_3$ &30&0.02s&113&2.43E-09&  &0.00s&25&0.00E+00\\
    spar030$\_100\_1$ &30&0.02s&114&2.36E-09& &0.01s&10&1.17E-13\\
    spar030$\_100\_2$ &30&0.02s&146&2.60E-09&  &0.00s&9&8.53E-14\\
    spar030$\_100\_3$ &30&0.01s&40&2.70E-09&  &0.00s&10&0.00E+00\\
    spar040$\_30\_1$ &40&0.04s&120&1.77E-09&  &0.01s&46&1.30E-09\\
    spar040$\_30\_2$ &40&0.02s&81&1.50E-09&  &0.01s&18&0.00E+00\\
    spar040$\_30\_3$ &40&0.01s&26&0.00E+00&  &0.00s&14&0.00E+00\\
    spar040$\_40\_1$ &40&0.01s&70&0.00E+00&  &0.00s&14&0.00E+00\\
    spar040$\_40\_2$ &40&0.03s&144&1.85E-09&  &0.00s&14&1.78E-15\\
    spar040$\_40\_3$ &40&0.02s&134&6.25E-10&  &0.00s&38&3.35E-14\\
    spar040$\_50\_1$ &40&0.06s&409&1.94E-09&  &0.00s&26&2.02E-13\\
    spar040$\_50\_2$ &40&0.02s&87&0.00E+00&  &0.00s&14&0.00E+00\\
    spar040$\_50\_3$ &40&0.03s&154&0.00E+00&  &0.00s&24&0.00E+00\\
    spar040$\_60\_1$ &40&0.05s&202&0.00E+00&  &0.01s&41&0.00E+00\\
    spar040$\_60\_2$ &40&0.03s&155&1.27E-09&  &0.00s&26&7.11E-14\\
    spar040$\_60\_3$ &40&0.02s&106&7.11E-09&  &0.00s&30&4.81E-11\\
    spar040$\_70\_1$ &40&0.03s&157&1.46E-09&  &0.01s&57&1.19E-09\\
    spar040$\_70\_2$ &40&0.01s&30&1.17E-09&  &0.00s&19&0.00E+00\\
    spar040$\_70\_3$ &40&0.06s&329&1.21E-09&  &0.01s&42&4.97E-14\\
    spar040$\_80\_1$ &40&0.01s&19&1.45E-09&  &0.00s&14&0.00E+00\\
    spar040$\_80\_2$ &40&0.02s&136&1.45E-09&  &0.00s&35&2.84E-14\\
    spar040$\_80\_3$ &40&0.02s&120&1.65E-09&  &0.00s&19&3.73E-14\\
    spar040$\_90\_1$ &40&0.03s&130&0.00E+00&  &0.00s&28&2.84E-14\\
    spar040$\_90\_2$ &40&0.03s&168&1.55E-09&  &0.00s&54&1.92E-13\\
    spar040$\_90\_3$ &40&0.04s&195&1.46E-09&  &0.01s&64&1.54E-09\\
    spar040$\_100\_1$ &40&0.03s&156&0.00E+00&  &0.00s&12&0.00E+00\\
    spar040$\_100\_2$ &40&0.02s&115&1.59E-09&  &0.00s&30&6.59E-11\\
    spar040$\_100\_3$ &40&0.06s&327&1.65E-09&  &0.00s&24&6.04E-14\\
    spar050$\_30\_1$ &40&0.01s&22&0.00E+00&  &0.00s&12&0.00E+00\\
    spar050$\_30\_2$ &40&0.02s&88&0.00E+00&  &0.00s&29&0.00E+00\\
    spar050$\_30\_3$ &40&0.02s&80&8.45E-10&  &0.00s&35&2.68E-10\\
    spar050$\_40\_1$ &40&0.05s&200&9.66E-10&  &0.01s&67&3.58E-10\\
    spar050$\_40\_2$ &40&0.02s&131&1.26E-09&  &0.00s&30&2.13E-14\\
    spar050$\_40\_3$ &40&0.02s&119&1.29E-09&  &0.00s&27&1.95E-14\\
    spar050$\_50\_1$ &50&0.05s&275&1.27E-09&  &0.01s&88&1.32E-09\\
    spar050$\_50\_2$ &50&0.03s&181&1.30E-09&  &0.01s&70&3.65E-10\\
    spar050$\_50\_3$ &50&0.06s&307&1.42E-09&  &0.00s&24&1.76E-13\\
    spar060$\_20\_1$ &60&0.02s&60&0.00E+00&  &0.00s&19&0.00E+00\\
    spar060$\_20\_2$ &60&0.06s&209&2.78E-13&  &0.01s&38&2.31E-10\\
    spar060$\_20\_3$ &60&0.06s&201&0.00E+00&  &0.00s&31&0.00E+00\\
    \bottomrule
    \end{tabular}

\end{table}
%%%%%%%%%%%%%%%%%%%%%%%%%%%%%%%%%%%%%%%%%%%%%%%%%%%%%%%%%%%%%%%%%%%%%%%%%%%%%
\begin{figure}[!htb]
  \centering
  \vspace{-3cm}
  \includegraphics[width=7cm]{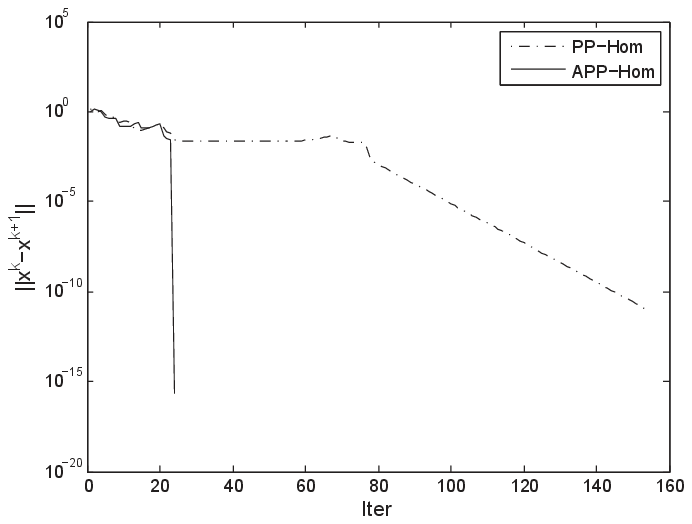}\hspace{-2.5cm}
  \includegraphics[width=7cm]{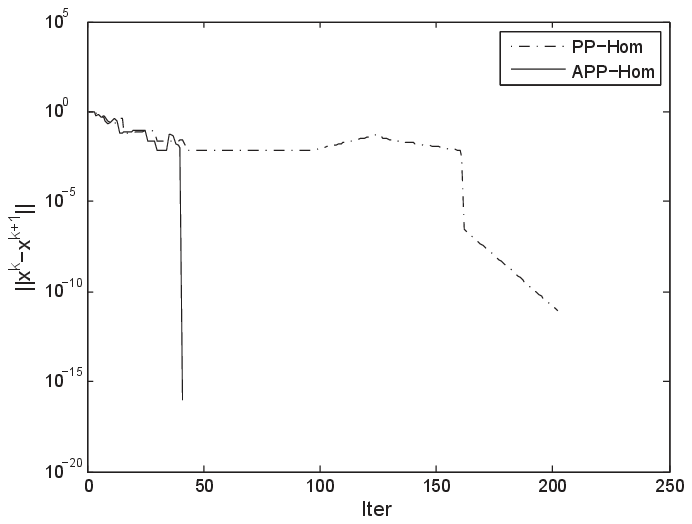}\\
  \vspace{-3cm}
  \caption{Iterative points distance for spar$\_$040$\_$050$\_$3 (left) and  spar$\_$040$\_$060$\_$1 (right).}\label{figure3}

\end{figure}

%%%%%%%%%%%%%%%%%%%%%%%%%%%%%%%%%%%%%%%%%%%%%%%%%%%%%%%%%%%%%%%%%%%%%%%%%%%%%%

We compared PP-Hom and APP-Hom algorithms with LBFGS-B, ASA, MINQ8 and TRR on solving NCBQP-D, NCBQP-S and  non-convex BQPs from \cite{vandenbussche2005branch}. The parameters in (\ref{switch condition}) is set as
$$\varepsilon=10^{-1}~ {\rm and}~a=2.$$
The results show APP-Hom outperforms LBFGS-B, ASA, MINQ8 and TRR. Fig. \ref{figure2} and \ref{figure3} shows that PP converges linearly at local and APP exhibits obvious acceleration than PP. Moreover, the minimize eigenvalue of $Q_{\bar{\mathcal{A}}_{m}\bar{\mathcal{A}}_{m}}$ is bigger than 0, which verifies Theorem \ref{theorem 2}.

$\bullet$  \textbf{SVM for recognition.} In this part, we tested the homotopy algorithm on solving BQPs in SVM \cite{sra2012optimization} for digit recognition and speech recognition. Given training set $\{X_i,y_i\}_{i=1}^{n}$ and testing set $\{T_j,s_j\}_{j=1}^{n}$, where $X_i, T_j$ are feature vectors and $y_i, s_j\in\{-1,+1\}$ are the labels. SVM classifies the testing set by a classifier
$f(x)={\rm sign}(\sum_{i=1}^{n}y_i\alpha_i^{*}K(X_i,x)+b^*),$
where $K$ is called kernel function, $b^*=y_j-\sum_{i=1}^{n}y_i\alpha_i^{*}K(X_i,X_j), \alpha_{j}^{*}>0$, and $\alpha^*$ is the solution of the following optimization:
\begin{eqnarray} \nonumber
&&\min_{\alpha} ~~~~\frac{1}{2}\sum_{i=1}^{n}\sum_{i=1}^{n}y_iy_j\alpha_i\alpha_j K(X_i,X_j)-\sum_{i=1}^{n}\alpha_i\\\label{svm}
&&s.t.~~~~\sum_{i=1}^{n}y_i\alpha_i=0,\\\nonumber
&& ~~~~~~~~~0\leq\alpha_i\leq C, i=1,...,n.
\end{eqnarray}

Based on the framework of augmented Lagrangian method (ALM), we solved (\ref{svm}) by solving every subproblem of ALM with the homotopy algorithm exactly and compared with the QP solver (interior-point method, denoted by IPM) of CPLEX 12.6 and the sequential minimal optimization (SMO) method \cite{fan2005working} in LIBSVM 3.22 \cite{CC01a}.

We did experiments with two database. The first database is the isolated letter speech database from UCI \cite{Lichman2013}, which contains training set with 6238 samples and testing set with 1559 samples. This database has 26 classifications: A-Z, and every sample has 617 attributes.  The second one is the  mnist database of handwritten digits\footnotemark[2]\footnotetext[2]{http://yann.lecun.com/exdb/mnist/}, which contains training set with 60000 samples and testing set with 10000 samples. Every sample is one $28\times 28$ pixel figure, that is, every sample has 784 attributes. This database has ten classifications as follows:
\begin{figure}[!htb]
  \centering
  \vspace{-0cm}
  \includegraphics[width=9.5cm]{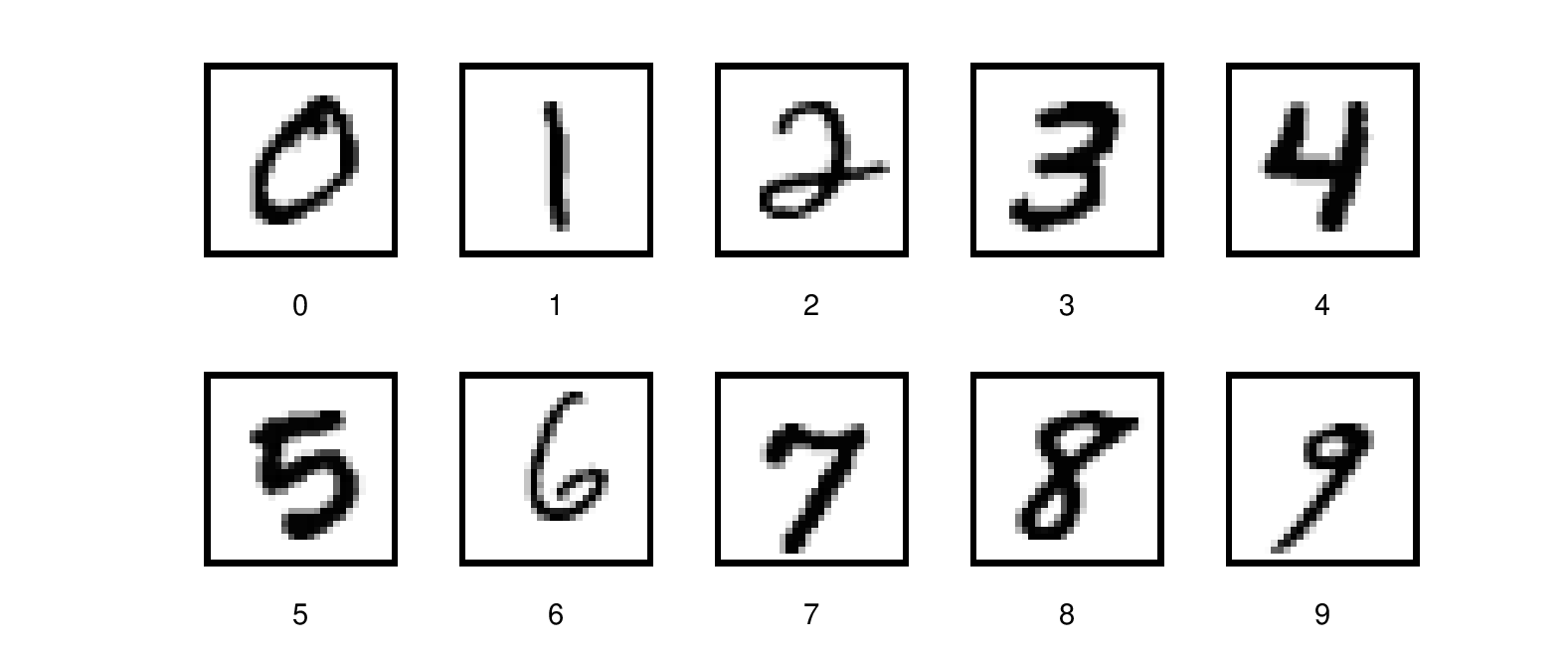}\\

  \caption{Ten classifications of the mnist handwritten digits database}
  \vspace{0cm}
\label{figure4}
\end{figure}

($\ref{svm}$) is a 2-class model and the letter speech and mnist databases are multi-class problems. So we did multi-classification problems with the following two strategies respectively.

The first strategy is that, for any classification $p, p=1,..,P$, where $P$ denotes the number of the classifications, we obtain $\alpha^{*,p}, b^{*,p}$ by solving ($\ref{svm}$) with
\begin{eqnarray} \nonumber
y^p=\left\{
\begin{array}{lcl}
1,      &      & X_i\in cl.p;\\
-1,  &      & else,
\end{array} \right.
\end{eqnarray}
then we give classifier as follows:
\begin{eqnarray} \label{classifier}
f_1(x)=\arg\max_{p}~ (y^p_i\alpha_i^{*,p}K(X_i,x)+b^{*,p}).
\end{eqnarray}

The second strategy is that, for any $p\neq q\in\{1,...,P\}$, choose the sample from the training sets which corresponding label is $p$ or $q$, and let
$$y^{p,q}=\left\{
\begin{array}{lcl}
1,      &      & X_i\in cl.p;\\
-1,  &      &  X_i\in cl.q,
\end{array} \right.$$
then we  give the classifier as follows:
\begin{eqnarray} \label{classifier2}
f_2(x)=\arg\max_{p}~ \sum_{q}{\rm sign}(y^{p,q}_i\alpha_i^{*,p,q}K(X_i,x)+b^{*,p,q}).
\end{eqnarray}

The first strategy needs to solve $(\ref{svm})$ for $P$ times, and every problem is in a size of the whole training set. While the second strategy needs to solve $(\ref{svm})$ for $\frac{P(P-1)}{2}$ times, and every time it only needs to solve a problem with the size equalling to number of classification $p$ and $q$. In LIBSVM 3.22, the multi-classifier adopts the second strategy.

In our experiments, we used the polynomial kernel
\begin{eqnarray} \nonumber
K(x,y)=((x\cdot y)+c)^d
\end{eqnarray}
with $c=0$, $d=2$ to classify the spoken letter and the mnist databases.

In Table \ref{table11}-\ref{table12}, we list the results of our algorithms (we use ALM-Hom to denote ALM with the subproblems solved exactly by the homotopy algorithm), IPM and SMO method. All of these methods did not use heuristic technique before solving, that is, we started them from zeros,  where ``Err.1'' and ``Err.2'' respectively denote the number of misclassification of classifier $(\ref{classifier})$ and classifier $(\ref{classifier2})$. Here, we just listed the running time of the first strategy, and did not list the time of the second strategy for it contains $\frac{P(P-1)}{2}$ parts. But we give the total time of the second strategy in the title of the tables. Moreover, we did not list the results of IPM for the Mnist database for it took much more time than the other two algorithms. The results show that ALM-Hom outperforms IPM and is competitive to SMO method on solving these two databases. Moreover,  ALM-Hom can obtain high accuracy output if it is required, while SMO method is hard to produce.

\begin{table}[!htb]
 \caption{\scriptsize{Experiments on SVM optimization for classification of the isolated letter speech database. The time in this table is the running time of the first strategy. In the second strategy, ALM-Hom took 5.17s, IPM took 107.08s and SMO took 13.37s in all.}}

    \label{table11}
    \tiny

    \begin{tabular}{p{1cm}p{0.6cm}p{0.6cm}p{0.6cm}p{0.2cm}p{0.6cm}p{0.6cm}p{0.6cm}p{0.2cm}p{0.6cm}p{0.6cm}p{0.6cm}}

    \toprule
    \multirow{2}{*}{Class.}& \multicolumn{3}{c}{ALM-Hom} && \multicolumn{3}{c}{IPM (CPLEX)}&& \multicolumn{3}{c}{SMO (LIBSVM)} \\
    \cline{2-4}\cline{6-8}\cline{10-12}
    &Time/s&Err.1&Err.2&&Time/s&Err.1&Err.2&&Time/s&Err.1&Err.2\\
    \hline
    cl.A &3.71&0&0&&240.80&0&0&&6.65&0&0\\
    cl.B &4.91&5&4&&242.58&5&4&&7.23&5&4\\
    cl.C &1.99&0&0&&259.33&0&0&&5.22&0&0\\
    cl.D &5.28&4&3&&233.34&4&3&&6.86&4&3\\
    cl.E &3.55&0&2&&255.32&0&2&&6.77&0&2\\
    cl.F &3.46&0&2&&232.10&0&2&&6.84&0&2\\
    cl.G &4.45&0&0&&260.34&0&0&&6.98&0&0\\
    cl.H &2.08&0&0&&220.50&0&0&&4.82&0&0\\
    cl.I &2.04&1&1&&223.03&1&1&&5.37&1&1\\
    cl.J &2.90&1&1&&273.58&1&1&&6.45&1&1\\
    cl.K &4.17&2&2&&218.09&2&2&&6.94&1&2\\
    cl.L &1.95&0&0&&200.01&0&0&&5.33&0&0\\
    cl.M &3.74&9&7&&252.48&9&6&&5.34&9&6\\
    cl.N &5.52&8&9&&230.56&8&9&&6.73&8&9\\
    cl.O &2.30&0&0&&208.38&0&0&&7.28&0&0\\
    cl.P &7.25&0&6&&235.25&0&5&&5.41&0&5\\
    cl.Q &1.99&4&0&&226.73&4&0&&7.78&4&0\\
    cl.R &1.68&0&0&&258.80&0&0&&5.22&0&0\\
    cl.S &1.72&3&3&&229.37&3&3&&4.94&3&3\\
    cl.T &5.10&3&5&&231.48&3&6&&7.30&3&6\\
    cl.U &2.00&2&2&&225.77&2&2&&5.88&2&2\\
    cl.V &6.10&5&6&&231.16&5&5&&7.02&5&5\\
    cl.W &2.96&0&0&&256.07&0&1&&6.57&0&1\\
    cl.X &1.88&0&0&&228.58&0&0&&5.02&0&0\\
    cl.Y &1.67&0&0&&220.33&0&0&&4.49&0&0\\
    cl.Z &2.41&4&3&&247.07&4&3&&5.91&4&3\\
    Total &86.87&51&56&&6138.86&51&55&&161.59&51&55\\

    \bottomrule
    \end{tabular}

\end{table}

\begin{table}[!htb]
 \caption{\scriptsize{Results of solving SVM optimization for classification of Mnist database. The time in this table is the running time of the first strategy. In the second strategy, ALM-Hom took 443.22s and SMO took 282.43s in all.}}

    \label{table12}
     \tiny

    \begin{tabular}{p{0.7cm}p{1.5cm}p{1cm}p{1cm}p{0.7cm}p{1.5cm}p{1cm}p{1cm}}
    \toprule
    \multirow{2}{*}{Class.}& \multicolumn{3}{c}{ALM-Hom} && \multicolumn{3}{c}{SMO (libsvm)} \\
    \cline{2-4}\cline{6-8}
    &Time/s&Err.1&Err.2&&Time/s&Err.1&Err.2\\
    \hline
    cl.0 &977.85&10&7 &&939.29&10&8\\
    cl.1 &729.54&10&7 &&540.36&10&8\\
    cl.2 &2160.36&22&24 &&2336.39&22&24\\
    cl.3 &2536.25&23& 23&&3501.45&23&25\\
    cl.4 &2107.25&17& 16&&1492.39&17&16\\
    cl.5 &2149.31&19&20 &&2397.18&19&19\\
    cl.6 &1281.68&17&18 &&1057.11&17&18\\
    cl.7 &2063.71&22&25 &&1926.60&22&26\\
    cl.8 &3234.19&24&22 &&4028.26&24&23\\
    cl.9 &2662.05&29& 30&&4111.90&29&28\\
    Total. &19907.2&193&192&&22230.8&193&195\\
    \bottomrule
    \end{tabular}

    \vspace{-0.2cm}
\end{table}

From (\ref{equation4.14a})-(\ref{equation4.14b}), we know that the homotopy algorithm needs to solve two linear systems with size $|\mathcal{A}_{m,i}|$ at each step. Moreover, from the results of the experiments, we found the number of the support vectors is small, that is, $\alpha^*$ is sparse, which implies the homotopy algorithm needs to solve small-scale linear systems. This good property of $\alpha^*$ makes ALM-Hom preform well on solving SVM optimizations.

$\bullet$  \textbf{Optimization in physics.} We also solved three optimization problems in physics: the obstacle problem, the elastic-plastic torsion problem, and the journal bearing problem in \cite{more1991solution}, which are formulated in a form

\begin{equation}\label{partial problem}
  \left\{
  \begin{aligned}
        &\min~q(v)=\frac{1}{2}\int_{\mathcal{D}}w\|\nabla v\|^{2}d\mathcal{D}-\int_{\mathcal{D}}fvd\mathcal{D} \\
        &~~{\rm s.t.}~v\in K,
  \end{aligned}
  \right.
\end{equation}
where $\mathcal{D}$ is an open set with a reasonably smooth boundary $\partial \mathcal{D}$, and $w,f\in H_{0}^{1}(\mathcal{D})$.

The obstacle problem A is a case of (\ref{partial problem}) with
\begin{equation}\label{obstacle problemA}
  \left\{
  \begin{aligned}
        &\mathcal{D}=(0,1)\times(0,1), \\
        & K=\{ H_{0}^{1}(\mathcal{D}):v_l\leq v\leq v_u~on~\mathcal{D}\},\\
        &v_l(x,y)=\sin(3.2x)\sin(3.2y), v_u(x,y)=2000;\\
        &w\equiv 1, f\equiv c,
  \end{aligned}
  \right.
\end{equation}
the obstacle problem B with
\begin{equation}\label{obstacle problemB}
  \left\{
  \begin{aligned}
        &\mathcal{D}=(0,1)\times(0,1), \\
        &K=\{ H_{0}^{1}(\mathcal{D}):v_l\leq v\leq v_u~on~\mathcal{D}\},\\
        &v_l(x,y)=(\sin(9.3x)\sin(9.3y))^3, \\
        &v_u(x,y)=(\sin(9.3x)\sin(9.3y))^2+0.02\\
        &w\equiv 1, f\equiv c,
  \end{aligned}
  \right.
\end{equation}
the elastic-plastic torsion problem with
\begin{equation}\label{elastic problem}
  \left\{
  \begin{aligned}
        &\mathcal{D}=(0,1)\times(0,1), \\
        &K=\{v\in H_0^1(\mathcal{D}):|v(x)|\leq {\rm{dist}}(x,\partial \mathcal{D}),~x\in \mathcal{D}\},\\
        &w\equiv 1, f\equiv c,
  \end{aligned}
  \right.
\end{equation}
and the journal bearing problem with
\begin{equation}\label{journal problem}
  \left\{
  \begin{aligned}
        &\mathcal{D}=\{(\theta, y):0<\theta<2\pi,0<y<2b\}, \\
        &K=\{v\in H_0^1(\mathcal{D}):v\geq 0~{\rm on}~\mathcal{D}\},\\
        &w=(1+\epsilon \cos\theta)^3, f=\epsilon \sin \theta, \epsilon \in (0,1).
  \end{aligned}
  \right.
\end{equation}

We followed Mor{\'e} \cite{more1991solution} by using finite difference to discretize (\ref{partial problem}). For convenience, assume
$\mathcal{D}=(d_1,d_2)\times(d_3,d_4)$. Let
$h_x$ and $h_y$ denote the grid spacings and
$$z_{i,j}=(d_1+ih_x,d_3+jh_y),~0\leq i\leq n_x+1,~0\leq j\leq n_y+1$$
denote the grid points, where
$n_x=\frac{d_2-d_1}{h_x}, n_y=\frac{d_4-d_3}{h_y}.$

Then we have
\begin{eqnarray}\nonumber
\int_{\mathcal{D}}w\|\nabla v\|^{2}d\mathcal{D}&=&\frac{h_x h_y}{4}\sum_{i,j}(\mu_{i,j}(\frac{v_{i+1,j}-v_{i,j}}{h_x})^2+\mu_{i,j}(\frac{v_{i,j+1}-v_{i,j}}{h_y})^2\\\nonumber
&+&\lambda_{i,j}(\frac{v_{i-1,j}-v_{i,j}}{h_x})^2+\lambda_{i,j}(\frac{v_{i,j-1}-v_{i,j}}{h_y})^2)
\end{eqnarray}
and
\begin{eqnarray}\nonumber
\int_{\mathcal{D}}fvd\mathcal{D}=h_x h_y\sum f_{i,j}v_{i,j},
\end{eqnarray}
where
$\mu_{i,j}=\frac{h_x h_y}{6}(w_{i+1,j}+w_{i,j}+w_{i,j+1})$ and $\lambda_{i,j}=\frac{h_x h_y}{6}(w_{i-1,j}+w_{i,j}+w_{i,j-1})$.

The above discrete problems have been included in the CUTEr test set\footnotemark[4] \footnotetext[4]{http://www.cuter.rl.ac.uk/Problems/mastsif.shtml} (e.g. OBSTCLAE, TORSION1, JNLBRNGA, etc). We solved the discrete problems by the homotopy algorithm and compared with LBFGS-B, ASA, MINQ8 and TRR. We use $n=n_x n_y$ to denote the dimension of the discrete problems. The results show that the homotopy algorithm took less time but obtained higher precision solutions than LBFGS-B, ASA, MINQ8 and TRR.

%%%%%%%%%%%%%%%%%%%%%%%%%%%%%%%%%%%%%%%%%%%%%%%%%%%%%%%%%%%%%%%%%%%%%%%%%%%%
\begin{table}[!htb]
\caption{\scriptsize{Five algorithms solving the discrete partial differential optimization problems.}}
\label{table13}
 \tiny
\begin{tabular}{cccccccccccccccc}
\toprule
\multirow {2}{*}{Problem}&
\multirow {2}{*}{n}&
\multicolumn{1}{c}{Homotopy} &
\multicolumn{1}{c}{LBFGS-B}&
\multicolumn{1}{c}{TRR} &
\multicolumn{1}{c}{ASA}  &
\multicolumn{1}{c}{MINQ8}  \\
&&Time $|$ $\| g(x)\|$ & Time$|$ $\| g(x)\|$ & Time$|$$\| g(x)\|$ & Time$|$$\| g(x)\|$& Time$|$$\| g(x)\|$\\
\midrule
\multirow {3}{*}{OA}
   &6400  &0.19s$|$1.50E-14&~~0.42s$|$1.79E-07&0.35s$|$3.26E-07 &0.35s$|$3.26E-08&18.33s$|$3.10E-07\\

   &10000 &0.22s$|$1.41E-14&~~0.78s$|$1.17E-07&0.90s$|$5.42E-07& 0.76s$|$4.11E-08&107.21s$|$4.23E-06 \\

   &14400 &0.88s$|$6.16E-14&~~2.90s$|$2.24E-07&3.44s$|$2.42E-06& 1.97s$|$4.23E-08&342.17s$|$7.19E-06 \\
\midrule
\multirow {3}{*}{OB}
   &6400  &0.78s$|$5.35E-14&~~4.23s$|$4.37E-02&1.89s$|$1.02E-07& 2.61s$|$8.38E-08&422.33s$|$4.28E-04\\

   &10000 &0.90s$|$3.66E-13&~~6.79s$|$2.83E-02& 3.55s$|$6.03E-08& 6.28s$|$8.92E-08&2236.11s$|$3.66E-04\\

   &14400 &1.55s$|$4.35E-13 &25.55s$|$4.16E-02 &8.64s$|$5.92E-07& 12.32s$|$1.17E-08&4572.53s$|$3.18E-04\\
\midrule
\multirow {3}{*}{EPT}
   &6400      &0.42s$|$4.17E-14&~~2.68s$|$8.58E-08&1.12s$|$1.07E-07& 0.48s$|$1.77E-08&216.04s$|$4.10E-05 \\

   &10000 &0.93s$|$1.45E-13&~~5.01s$|$2.81E-06&2.38s$|$1.79E-07& 0.89s$|$6.95E-08&1212.44s$|$7.32E-05\\

   &14400  &1.49s$|$3.44E-11&18.59s$|$1.29E-05&4.06s$|$1.54E-07& 1.82s$|$1.20E-08&2429.37s$|$8.82E-05 \\
\midrule
\multirow {3}{*}{JB}
   &6400  &0.32s$|$8.07E-14&~~2.82s$|$5.80E-07&0.51s$|$1.62E-06& 0.67s$|$9.12E-08&13.67s$|$2.42E-05\\

   &10000 &0.54s$|$1.36E-13&~~4.44s$|$1.05E-05&1.14s$|$2.54E-06& 1.22s$|$1.10E-07&68.91s$|$3.91E-05\\

   &14400  &1.04s$|$1.76E-13&16.42s$|$3.33E-05&2.82s$|$2.32E-06 & 2.04s$|$1.09E-07&121.77s$|$7.22E-05\\

   \bottomrule
\end{tabular}

\vspace{-0.2cm}
\end{table}

\begin{figure}[!htb]
  \centering
  \vspace{-0cm}
  \includegraphics[width=5cm]{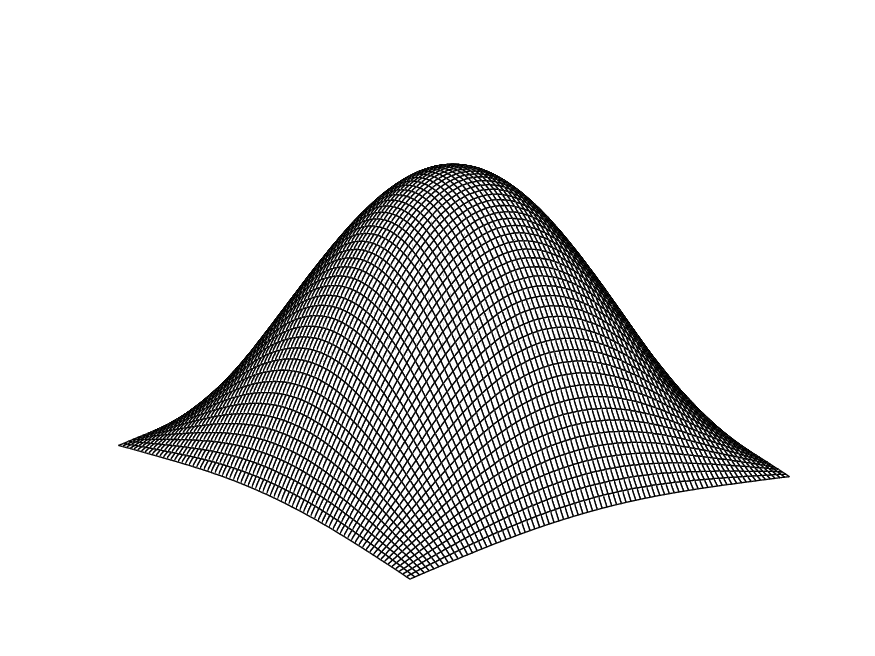}
  \hspace{-0cm}
  \includegraphics[width=5cm]{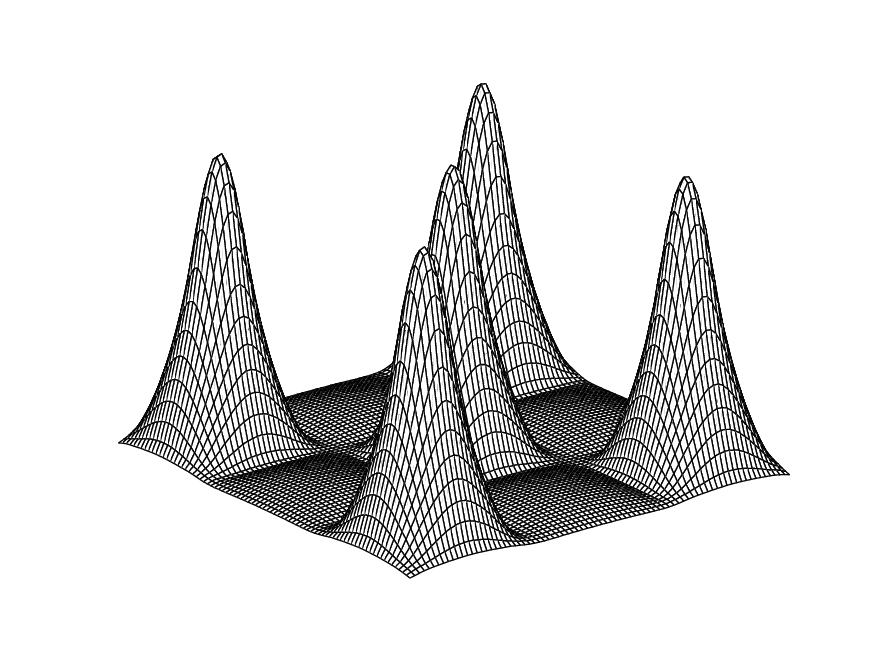}\\
  \vspace{-0.5cm}
  \includegraphics[width=5cm]{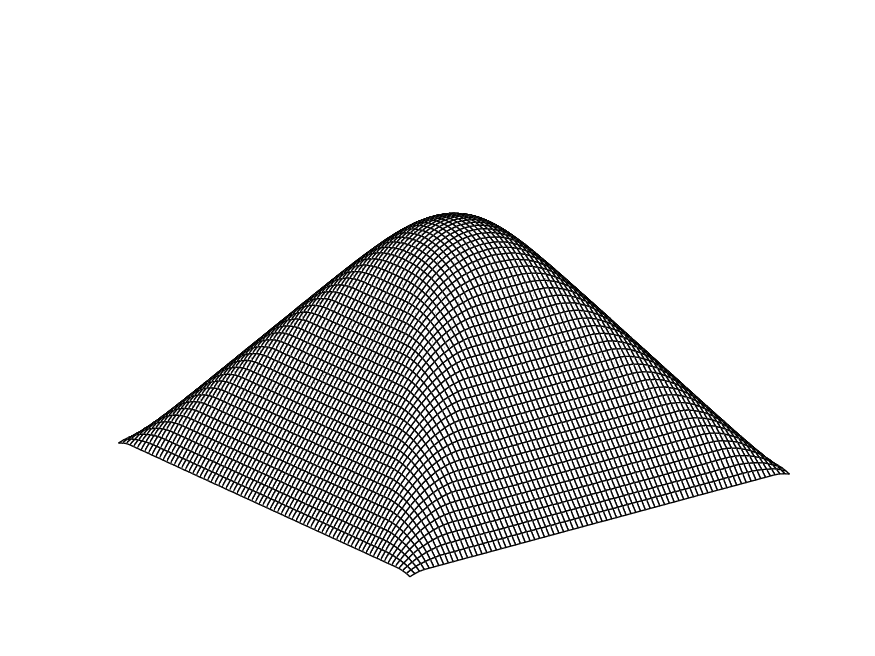}
  \hspace{-0cm}
  \includegraphics[width=5cm]{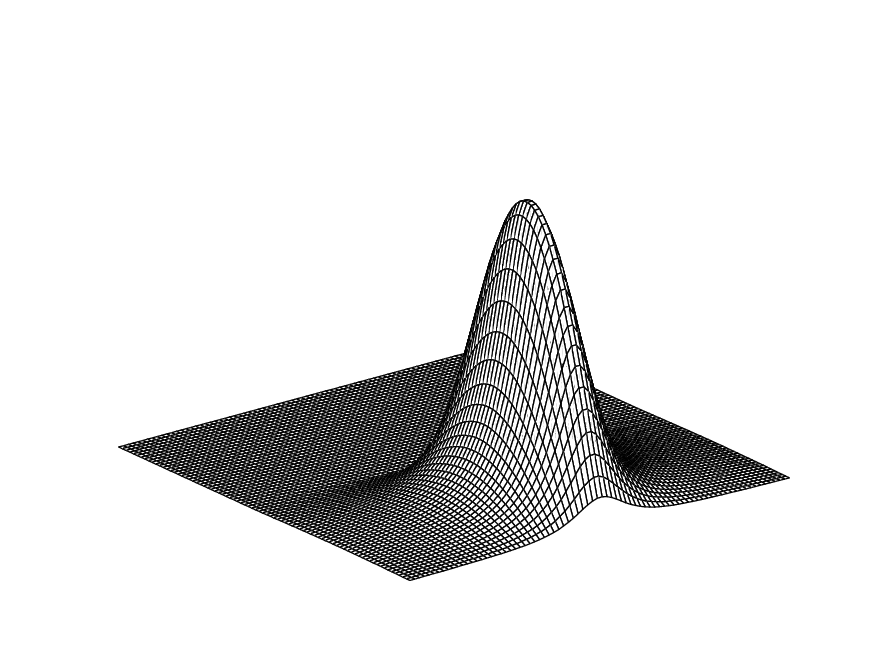}\\

  \caption{Obstacle problem A with c=1 (upper left), obstacle problem B with c=5 (upper right), elastic-plastic torsion problem with c=10 (lower left) and journal bearing problem with $\varepsilon$=0.8 (lower right).}
  \vspace{-0cm}
\label{figure5}
\end{figure}

Finally, we compared the homotopy algorithm with the PAS algorithm in qpOASES on solving the discrete PDE optimization problems from the approximate solution $\hat{z}$ and display the results in Table \ref{table14}. The results show that the homotopy method is much faster than the  PAS method in qpOASES.

%%%%%%%%%%%%%%%%%%%%%%%%%%%%%%%%%%%%%%%%%%%%%%%%%%%%%%%%%%%%%%%%%%%%%%%%%%%%%%%
\begin{table}[!htb]
\caption{\scriptsize{The homotopy algorithm and the PAS method in qpOASES solving the discrete partial differential optimization problems from the initial solution $\hat{z}$ like (\ref{equation4.24})}}
\label{table14}

\tiny
\begin{tabular}{p{2.8cm}p{1.6cm}p{1.2cm}p{1.2cm}p{1.2cm}p{1.2cm}p{1.2cm}}
\toprule
\multirow {2}{*}{Problem}&
\multirow {2}{*}{n}&
\multicolumn{2}{c}{PAS(qpOASES)} &
\multicolumn{2}{c}{Homotopy}

  \\

   &&Iter &Time & Iter& Time  \\
\midrule
\multirow {3}{*}{Obstacle A}
   &6400      &34&3.69s&29&0.10s \\

   &10000 &26&4.69s&22&0.14s \\

   &14400  &69&29.13s&66&0.54s \\
\midrule
\multirow {3}{*}{Obstacle B}
   &6400      &64&8.18&64&0.36s \\

   &10000 &90&23.67&88& 0.50s\\

   &14400  &97&43.21&89& 0.86s\\
\midrule
\multirow {3}{*}{Elastic-plastic torsion}
   &6400      &97&5.43&94&0.21 \\

   &10000 &103&14.98&96&0.54 \\

   &14400  &111&49.34&102&0.84 \\
\midrule
\multirow {3}{*}{Journal bearing}
   &6400      &49&2.33&41&0.14 \\

   &10000 &77&14.67&71& 0.29\\

   &14400  &93&34.45&87&0.56 \\

   \bottomrule
\end{tabular}

\vspace{-0.2cm}
\end{table}

\section{Conclusion}
PP converges $R$-linearly for BQP, furthermore, it converges $Q$-linearly if the limit point satisfies the strict complementary conditions. More precisely, PP is  linear iteration when the free variables do not change. According to this property, an accelerated PP algorithm  which exhibits the obvious effect of acceleration than PP algorithm is presented. The accelerated PP can automatically identify the endgame, that is, when the optimal active set is founded, it converges to the limit point very quickly.

The performance of the homotopy algorithm depends on the initial solution. Then the approximately solving is indispensable. APG algorithm is a good method to predict the optimal active set, because it costs small at each step, converges with a rate $O(\frac{1}{k^2})$ and is easy to implement. However, APG converges slowly at end of the iterations, which hinders it to be an independent algorithm for high-precision solutions, so we terminate it when some criteria are satisfied. With the approximately solving stage, the steps of the homotopy algorithm is greatly reduced.  Moreover, the $\varepsilon$-precision verification and correction steps ensure the stability of the homotopy algorithm.

Finally, the numerical results confirm our theoretical results well and show APP and the homotopy algorithm are effective in practice. Moreover, from the numerical results of solving SVM optimizations, we see the homotopy algorithm can utilize the sparsity of the solution well, such as for SVM optimization, the homotopy algorithm only needs to solve linear systems with the a size close to the number of the support vectors.
%%%%%%%%%%%%%%%%%%%%%%%%%%%%%%%%%%%%%%%%%%%%%%%%%%%%%%%%%%%%%%%%%%%%%%%

%%%%%%%%%%%%%%%%%%%%%%%%%%%%%%%%%%%%%%%%%%%%%%%%%%%%%%%%%%%%%%%%%%%%%%%%%%%%%%%%%%

\section*{Acknowledgments}

The authors would like to thank the colleagues for their valuable suggestions that led to improvement in this paper. This research was supported by the National Natural Science Foundation of China (11571061, 11401075, 11701065
), and the Fundamental Research Funds for the Central Universities (DUT16LK05)

%\begin{acknowledgements}
%\end{acknowledgements}

%\bibliography{qpboxreferences}

\end{document}